\documentclass[final,onefignum,onetabnum]{siamart250211}
\usepackage{marvosym}
\usepackage[symbol]{footmisc}
\usepackage{ulem,url}
\usepackage{pgfplots,pgf}
\usetikzlibrary{pgfplots.groupplots}
\usepackage{graphicx}
\usepackage{wrapfig}
\usepackage{verbatim}
\usepackage[justification=centering]{subcaption}
\pgfplotsset{compat=newest}
\usepackage{amssymb}
\usepackage{float}
\usepackage{verbatim}
\usepackage[top=20mm, bottom=20mm, left=35mm, right=35mm]{geometry}
\usepackage{color}
\usepackage{bm}
\usepackage{algorithm}
\usepackage{algorithmicx}
\usepackage{algpseudocode}
\usepackage{booktabs}

\newcommand{\pvct}[1]{\bm{#1}}
\newcommand{\vct}[1]{\bm{\mathsf{#1}}}
\newcommand{\mtx}[1]{\bm{\mathsf{#1}}}

\newtheorem{thm}{Theorem}
\newtheorem{assumption}[thm]{Assumption}

\theoremstyle{definition}
\newtheorem{remark}{Remark}

\definecolor{MyDarkGreen}{rgb}{0,0.6,0}
\definecolor{MyGray}{rgb}{0.6,0.6,0.6}

\newcommand{\pxx}{\pvct{x}}

\newcommand{\RN}[1]
    {\MakeUppercase{\romannumeral #1}}

\begin{document}

\begin{center}
\textsc{An overlapping domain decomposition method based on\\ solution-transfer operators}

\vspace{3mm}
\renewcommand*{\thefootnote}{\fnsymbol{footnote}}
Simon Dirckx\footnote{Oden Institute, University of Texas at Austin. \Letter\;\texttt{simon.dirckx@austin.utexas.edu}}, Anna Yesypenko\footnote{Department of Mathematics, The Ohio State University. \Letter\;\texttt{yesypenko.1@osu.edu}}, Per-Gunnar Martinsson\footnote{Oden Institute, University of Texas at Austin. \Letter\;\texttt{pgm@oden.utexas.edu}}
\vspace{10mm}

\begin{minipage}{120mm}
\textsc{Abstract:}
An overlapping domain decomposition method is described for variable-coefficient elliptic boundary value problems on domains that can be decomposed into slabs or shells. The method represents the global solution through its traces on internal interfaces, coupled by local Dirichlet solution transfer operators posed on overlapping double slab domains. The key observation is that these interface maps act between separated interfaces, and as such can be written as smooth-kernel integral operators. The resulting global equilibrium system is Fredholm second kind and, unlike non-overlapping formulations, requires no same-interface Dirichlet-to-Neumann or other interface maps with singular kernels. This makes its off-diagonal blocks highly amenable to hierarchical low-rank compression. The formulation admits complementary continuum and discrete interpretations. At fixed slab width, the method is stable under discretization, sufficiently accurate local solves and compression. At the discrete level, the system can be interpreted as a block Jacobi preconditioned Schur complement system. 
For compatible SPD discretizations satisfying a standard stable-splitting assumption, the symmetrically scaled interface matrix satisfies an energy-norm condition-number bound that depends on the slab width but is uniform with respect to the local resolution.
The formulation is implemented using high-order local solvers and hierarchical compression based on randomized sampling. 
Numerical experiments report iteration counts, accuracy, and compressibility for 2D and 3D elliptic, nonsymmetric, and oscillatory problems with as many as 28 million degrees of freedom.
\end{minipage}
\end{center}
\renewcommand*{\thefootnote}{\arabic{footnote}}
\setcounter{footnote}{0}

\vspace{10mm}


\section{Introduction}
\begin{figure}[H]
\begin{center}
\includegraphics[width=.3\linewidth,height = 3cm]{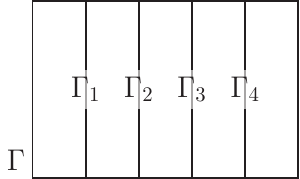}
\hfill
\includegraphics[width=.25\linewidth]{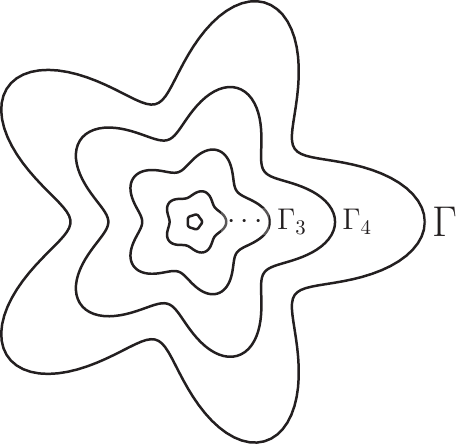}
\hfill
\includegraphics[width=.3\linewidth,height = 2.5cm]{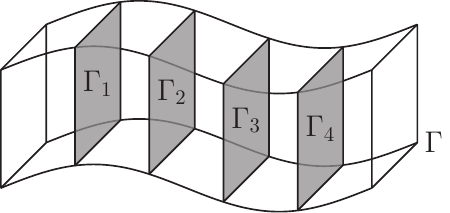}
\end{center}
\caption{Examples of domains that can naturally be tessellated into thin slabs or shells. Each is decomposed into $4$ overlapping double-wide-slabs, derived from $5$ non-overlapping slabs.}
\label{fig:geometries}
\end{figure}

We describe a numerical method for solving boundary value problems of the form
\begin{equation}
\label{eq:basic}
\left\{\begin{aligned}\mbox{}
[\mathcal{A}u](\pxx) =&\ f(\pxx)\,\qquad \pxx \in \Omega,\\
   u(\pxx) =&\ g(\pxx)\,\qquad \pxx \in \Gamma,
\end{aligned}\right.
\end{equation}
where $\Omega$ is a domain in $\mathbb{R}^{2}$ or $\mathbb{R}^{3}$
with boundary $\Gamma$, and where $\mathcal{A}$ is a scalar elliptic
partial differential operator that may have variable coefficients. We restrict ourselves in this manuscript to coefficients in $\mathbb{R}$, but generalizing to $\mathbb{C}$ poses no significant additional challenge.
The method described works best on domains that can naturally be tessellated into thin slabs, as illustrated in Figure \ref{fig:geometries}. The central object is our overlapping solution-transfer formulation, not a particular choice of local discretization. The solution-transfer blocks can in principle be approximated using finite differences, finite elements, spectral elements, or other local solvers, provided that the interface discretizations give accurate approximations of the local Dirichlet solution maps. For the large-scale numerical experiments in this paper, we use a multidomain spectral collocation implementation, since the thin local Dirichlet problems are particularly well suited to high-order discretizations and direct local solvers.

\subsection{A model problem}
\label{sec:toyproblem}
To introduce the key ideas, let us consider a model
problem where a square domain $\Omega$ has been partitioned into five
thin strips, as shown in Figure \ref{fig:FD}. We discretize (\ref{eq:basic})
to obtain a linear system $\mtx{A}\vct{u} = \vct{b}$ for some sparse
matrix $\mtx{A}$. For simplicity, assume that we use finite differences,
so that each entry of $\vct{u}$ holds a collocated
value of the solution $u$ at some grid point.
Now suppose that we --- in principle ---
use block Gaussian elimination to excise all nodes that are interior to
the five slabs (blue dots in Figure~\ref{fig:FD}), keeping only nodes associated with the interfaces.
This would result in a block tridiagonal linear system of the form
\begin{equation}
\label{eq:modelT}
\left[\begin{array}{rrrr}
\mtx{T}_{1,1} & \mtx{T}_{1,2} & \mtx{0}       & \mtx{0}       \\
\mtx{T}_{2,1} & \mtx{T}_{2,2} & \mtx{T}_{2,3} & \mtx{0}       \\
\mtx{0}       & \mtx{T}_{3,2} & \mtx{T}_{3,3} & \mtx{T}_{3,4} \\
\mtx{0}       & \mtx{0}       & \mtx{T}_{4,3} & \mtx{T}_{4,4}
\end{array}\right]
\left[\begin{array}{rrrr}
\vct{u}_{1} \\ \vct{u}_{2} \\ \vct{u}_{3} \\ \vct{u}_{4}
\end{array}\right]
=
\left[\begin{array}{rrrr}
\vct{f}_{1}' \\ \vct{f}_{2}' \\ \vct{f}_{3}' \\ \vct{f}_{4}'
\end{array}\right],
\end{equation}
where the blocks $\mtx{T}_{j,j'}$ are dense matrices formed by taking
Schur complements in the original sparse matrix $\mtx{A}$, where
the vectors $\vct{h}_{j}$ encode the body load and the boundary data,
and where the vectors $\vct{u}_{j}$ hold the function values of $u$ on
the grid nodes on the four internal boundaries (see Section \ref{sec:LA} for details).
Following standard practice, we could solve (\ref{eq:modelT}) using a
preconditioned iterative solver, where each block $\mtx{T}_{j,j'}$ is
applied implicitly using local solvers for the five domain interiors
(oftentimes a sparse direct solver, to accelerate repeated solves).

In the method proposed here, we consider the block diagonally pre-conditioned
version of the system (\ref{eq:modelT}), which results in a system of the form
\begin{equation}
\label{eq:modelS}
\left[\begin{array}{rrrr}
 \mtx{I}       & -\mtx{S}_{1,2} &  \mtx{0}       &  \mtx{0}       \\
-\mtx{S}_{2,1} &  \mtx{I}       & -\mtx{S}_{2,3} &  \mtx{0}       \\
 \mtx{0}       & -\mtx{S}_{3,2} &  \mtx{I}       & -\mtx{S}_{3,4} \\
 \mtx{0}       &  \mtx{0}       & -\mtx{S}_{4,3} &  \mtx{I}
\end{array}\right]
\left[\begin{array}{rrrr}
\vct{u}_{1} \\ \vct{u}_{2} \\ \vct{u}_{3} \\ \vct{u}_{4}
\end{array}\right]
=
\left[\begin{array}{rrrr}
\widehat{\vct{f}}_{1} \\ \widehat{\vct{f}}_{2} \\ \widehat{\vct{f}}_{3} \\ \widehat{\vct{f}}_{4}
\end{array}\right],
\end{equation}
where $\mtx{S}_{j,j'} = -\mtx{T}_{j,j}^{-1}\mtx{T}_{j,j'}$ and
$\widehat{\vct{f}}_{j} = \mtx{T}_{j,j}^{-1}\vct{f}_{j}'$. In contrast to standard practice, we will form the blocks $\mtx{S}_{j,j'}$ \textit{explicitly}, exploiting that they
have internal structure that allows us to store them to
high accuracy using data sparse representations
(see, e.g., \cite{hackbusch,2008_bebendorf_book,2002_hackbusch_H2,2004_borm_hackbusch,2010_borm_book,2005_martinsson_fastdirect,2012_martinsson_FDS_survey,2009_xia_superfast,2013_xia_randomized,2010_gu_xia_HSS,2016_kressner_rankstructured_review})
that exploit rank deficiencies in the off-diagonal blocks of $\mtx{S}_{j,j'}$. Importantly, we will not compute the blocks $\mtx{S}_{j,j'}$ using the formula $\mtx{S}_{j,j'} = -\mtx{T}_{j,j}^{-1}\mtx{T}_{j,j'}$, but we will use the fact that they can be written as Dirichlet-to-Dirichlet maps (see Section~\ref{sec:ConstructS}).

The key observation underpinning our work is that the
matrices $\mtx{S}_{j,j'}$ are very benign. They turn out
to be discrete approximations to integral operators that
are not only compact, but in fact have \textit{smooth}
kernels, with no singularity at the diagonal. In consequence, these matrices are highly
compressible. In practice, the matrix $\mtx{S}$ from the global interface system~\eqref{eq:modelS} is often substantially better conditioned than the original discretized PDE system. Its favorable behavior is connected to the Fredholm second-kind structure of the underlying continuum equilibrium equation. In discrete form, we write
\begin{equation}
\label{eq:pseudofredholm}
\mtx{S}\vct{u}=\bigl(\mtx{I} - \mtx{K}\bigr)\vct{u} = \widehat{\vct{f}},
\end{equation}
Here $\mtx{K}$ approximates a compact operator assembled from smooth-kernel maps between separated interfaces. This second-kind structure explains the observed favorable conditioning and, more precisely, why local refinement at fixed physical slab width does not introduce the large spectral scales associated with discretizing the underlying differential operator; Section~\ref{sec:continuum_trace} makes the latter approximation statement precise. 
For compatible SPD discretizations, Section~\ref{sec:conditioning} gives an upper bound for the symmetrically scaled system, equivalently for $\mtx{S}$ in an associated energy norm. The supplementary material reports empirical spectral diagnostics for representative stencil and collocation realizations of the raw matrix. For reasons that will become apparent in Section~\ref{sec:construct_reduced_system}, the matrix $\mtx{S}$ is referred to as the (discretized) \textit{equilibrium operator}.

\subsection{Constructing the reduced system}\label{sec:ConstructS}

A key observation of our work is that it is possible to directly
form the coefficient matrix in the linear system (\ref{eq:modelS}),
\textit{without first forming the blocks $\mtx{T}_{j,j'}$ in (\ref{eq:modelT})}.
We provide the details of this technique in Section \ref{sec:construct_reduced_system}, but the
idea is quite simple: the $j$'th block row  in (\ref{eq:modelS})
can be written as
\begin{equation}
\label{eq:localS}
\vct{u}_{j} = \mtx{S}_{j,j-1}\vct{u}_{j-1} + \mtx{S}_{j,j+1}\vct{u}_{j+1} + \widehat{\vct{f}}_{j},
\qquad j \in \{2,3\}.
\end{equation}
We refer to the matrices $\mtx{S}_{j,j-1}$ and $\mtx{S}_{j,j+1}$ as local solution-transfer operators. In the present Dirichlet-trace formulation, they are Dirichlet-to-Dirichlet operators that map trace data on $\Gamma_{j-1}$ and $\Gamma_{j+1}$ to solution values on the separated interface $\Gamma_j$. This means that we can consider a local Dirichlet problem on the subdomain $\Psi_{j}$ located between $\Gamma_{j-1}$ and $\Gamma_{j+1}$, cf.~Figure \ref{fig:local_bvp}, and construct $\mtx{S}_{j,j-1}$ and $\mtx{S}_{j,j+1}$ directly by simply solving this local problem.

In fact, the matrices $\mtx{S}_{j,j-1}$ and $\mtx{S}_{j,j+1}$ can be built using
\textit{any} local BVP solver defined on the double-wide slab $\Psi_{j}$.  The only requirement is that the discretizations of the local maps are sufficiently accurate.  In our implementation, available at \url{https://github.com/SimonDirckx/SslabLU}, we employ a multidomain spectral collocation method. Since the local domains are thin strips, sparse direct solvers are highly efficient even for problems in three dimensions.

In principle, it is not necessary to explicitly form the blocks $\mtx{S}_{j,j'}$. If an LU-factorization of the local stiffness matrix $\mtx{A}_{j}$ is computed, this can be used to compute the action of $\mtx{S}_{j,j+1}$ and $\mtx{S}_{j,j-1}$. We will refer to this approach as the ``\textit{matrix-free}'' approach. However, there are a number of scenarios in which one may wish to form the off-diagonal blocks $\mtx{S}_{j,j'}$ explicitly:
\begin{enumerate}
    \item the number of matrix vector-products with $\mtx{S}$ is very high (e.g., $\mtx{S}$ is used in an iterative solver with many right-hand sides and/or iterations),
    \item the memory footprint of storing all local factorizations is prohibitively large, or
    \item a direct solver for the equilibrium system is to be used.
\end{enumerate}
We return to the first two scenarios in Section~\ref{subsec:compare_to_matfree}. The development of a direct solver for $\mtx{S}$ is currently work in progress. 

In case the off-diagonal blocks are formed explicitly, the final component that enables high computational efficiency even in 3D is that, while the blocks $\mtx{S}_{j,j'}$ are all dense, they can be represented efficiently by exploiting that their off-diagonal blocks have low numerical rank. In this work, we use the Hierarchically Block Separable format of \cite{2012_martinsson_FDS_survey}. This format admits fast and simple matrix arithmetic, but its weak-admissibility structure is generally not effective for three-dimensional same-interface boundary operators. Same-interface Dirichlet-to-Neumann and related maps retain a local pseudodifferential singularity. They remain compressible in hierarchical formats with suitable admissibility conditions, but their near-field and adjacent-cluster interactions do not generally have uniformly low rank. By contrast, the matrices $\mtx{S}_{j,j'}$ act between separated interfaces and approximate integral operators whose kernels are smooth on the entire source-target product. This permits effective use of weak-admissibility HBS even in three dimensions. We obtain the rank
structured representations of $\mtx{S}_{j,j'}$ using the black-box randomized
compression technique of \cite{fastHBS}.

\begin{remark}
The use of rank-structured matrices in the context of problems with oscillatory
solutions is a delicate matter. It is well known that as the wave-length shrinks, or the solid angle between domains increases relative to the wave number, the ranks of the off-diagonal blocks increase (see~\cite{Dirckx}). This eventually makes $\mathcal{H}$-matrix arithmetic infeasible. For fixed slab width measured in wavelengths, the required ranks remain controlled in two dimensions (even to the point that the \textit{exact} rank often is small, see~\cite[Ch.~21]{2019_martinsson_fast_direct_solvers}). In three dimensions, the experiments show mild rank growth over the tested range.
\end{remark}

\subsection{Main features of the reduced formulation}
The reduced formulation has several useful features:

\textit{Fredholm second-kind structure.}
For fixed slab width $H$, the equilibrium operator has the form identity minus compact kernel integral operators between separated interfaces. As such, its kernels are smooth everywhere at the continuum level, and the resulting operator has the character of a second-kind trace equation.

\textit{Conditioning in the compatible SPD setting.}
For compatible SPD discretizations, the equilibrium matrix $\mtx{S}$ is self-adjoint and positive definite in the associated $\mtx{D}$-energy inner product. Equivalently, the symmetrically scaled matrix $\widetilde{\mtx{S}}=\mtx{D}^{-1/2}\mtx{T}\mtx{D}^{-1/2}$ satisfies $\kappa_2(\widetilde{\mtx{S}})\le 2C_{\rm ss}H^{-2}$. The scaling $\mtx{D}$ and its induced norm depend on the compatible discretization, but the bound is uniform with respect to local refinement at fixed $H$. The supplementary material reports empirical diagnostics suggesting that the raw Euclidean condition number follows a similar trend in the stencil and collocation realizations tested. Nonsymmetric and oscillatory problems are studied empirically.

\textit{Data sparsity.}
The off-diagonal blocks $\mtx{S}_{j,j \pm 1}$ in the reduced system $\mtx{S}$ are discrete Dirichlet-to-Dirichlet operators between separated interfaces.  This is the reason that formats based on weak admissibility, such as HODLR or HBS/HSS, are effective even in three dimensions.

\textit{Explicit compressed assembly.}
A matrix-free implementation would apply each off-diagonal block $\mtx{S}_{j,j \pm 1}$ of the equilibrium system $\mtx{S}$ by performing local solves at every Krylov iteration. We instead construct compressed representations of the dense blocks during the setup phase. This pays the local-solve cost up front, enables fast repeated matrix-vector products, and is particularly attractive for multiple right-hand sides or for a future direct solver for the equilibrium system.

\textit{High-order local discretizations and parallelism.}
The local problems are posed on thin slabs, so sparse direct solvers are effective even when high-order (multidomain) spectral discretizations are used. The local LU factorizations and blocks compression are independent from slab to slab, and matrix-vector products with the assembled operator are naturally parallel.

\subsection{Connection to prior work and classical domain decompositions}
The method outlined in this manuscript belongs to the broad class of \textit{substructuring methods}; a particular type of domain decomposition methods, in which interior degrees of freedom are eliminated and a reduced problem is posed as an equilibrium equation using interface maps. Standard references include \cite{smithDD,toselliWidlundDD}. Classical substructuring methods arise from non-overlapping formulations where the interface maps are discrete Poincar\'e--Steklov maps. These can be Dirichlet-to-Neumann maps, as in Section~\ref{sec:toyproblem}, but Robin-to-Robin and impedance-to-impedance operators can also be used. For classical substructuring methods, the resulting equilibrium systems contain self-interactions, i.e., maps between interfaces and themselves. 
Because these operators are pseudodifferential, their same-interface kernels are locally singular and are less amenable to weak-admissibility compression.
Appendix~B compares the Dirichlet-to-Neumann, impedance-to-impedance, and overlapping equilibrium formulations for the numerical example introduced in Section~\ref{sub:example_2_photonic_crystal_waveguide}. In contrast to the classical non-overlapping substructuring methods, our method considers an overlapping formulation that only ever constructs interface maps between separated interfaces. 
In this sense, it is closely related to the local Dirichlet solves underlying Schwarz' original method~\cite{Schwarz1870}.
Recently, there has been renewed interest in overlapping substructuring methods, as in, e.g.,~\cite{Ciaramella2022}. 
In contrast to the method there, we consider a broad class of local discretizations and study the compressibility of the resulting interface operators.
We do not consider Schwarz iterations or coarse space correction here. In contrast, we iterate on the global equilibrium system.

The  discrete equilibrium system $\mtx{S}$ can be viewed algebraically as a preconditioned Schur-complement, but its blocks $\mtx{S}_{j,j\pm1}$ also have a direct continuum interpretation as local Dirichlet-to-Dirichlet maps: solution operators mapping data from one interface to solution values on a \textbf{separated} interface. This separation removes the kernel singularity present in the non-overlapping formulation.

For wave problems, domain decomposition methods often use transmission conditions designed to reduce internal numerical resonances. This includes optimized Schwarz methods, sweeping preconditioners, polarized traces and impedance-to-impedance maps \cite{sweeping_hmat,sweeping_pml,ZEPEDANUNEZ2016347,Zepeda-Nunez2019-lp,doi:10.1137/16M109781X,2013_martinsson_ItI,Despres1991}. 
For oscillatory problems, the solution-transfer formulation separates refinement-induced stiffness from frequency-dependent stability. At fixed frequency and fixed slab width, local refinement approximates a fixed interface equation built from local solution maps. Frequency-dependent effects, including proximity to local or global resonances, remain. The Helmholtz experiments investigate this behavior empirically. Combining the separated-interface viewpoint with wave-adapted trace variables is a natural direction for future work.

The closest precursor to the present work is ``SlabLU'', the slab-based direct solver of \cite{Yesypenko_2024} which also uses high-order multidomain solvers, randomized compression, and rank-structured $\mathcal{H}$-matrix algebra, but uses a non-overlapping slab decomposition. 
The present work introduces overlap and replaces the same-interface blocks of the non-overlapping formulation by Dirichlet-to-Dirichlet solution-transfer maps between separated interfaces.

\subsection{Outline}
Our manuscript is structured as follows: in Section~\ref{sec:construct_reduced_system} the reduced equilibrium system is derived numerically and analytically. Its construction using local solvers and randomized compression is outlined. Section~\ref{sec:continuum_trace} establishes the continuum second-kind structure and stability under convergent trace approximation.
Section~\ref{sec:conditioning} analyzes compatible SPD discretizations and proves a conditioning bound in the natural energy norm that depends on the slab width, but not on the local resolution.
In section~\ref{sec:complexity} we analyze the complexity of our method and compare it to the non-overlapping and matrix-free approach. In section~\ref{sec:numerical} we test our method on illustrative and challenging numerical experiments.

\section{Derivation of the reduced linear system}

\label{sec:construct_reduced_system}

This section describes how we construct the reduced linear system (\ref{eq:modelS})
that forms the foundation of our method.
Throughout this section, we assume that the computational domain $\Omega$ in (\ref{eq:basic})
has been tessellated into $N_{\text{ds}}+1$ thin slabs separated by some interfaces $\{\Gamma_{j}\}_{j=1}^{N_{\text{ds}}}$, forming $N_{\text{ds}}$ double slabs $\{\Psi_j\}_{j=1}^{N_{\text{ds}}}$ (cf.~Figure \ref{fig:geometries}).
Our objective is to build the blocks in the coefficient matrix in (\ref{eq:modelS}),
as well as the equivalent reduced loads. We restrict our attention in this manuscript to the case of equispaced slabs. For irregularly spaced slabs the derivations in this Section as well as Section~\ref{sec:continuum_trace} carry through unchanged, but we make no claim about the conditioning (see Section~\ref{sec:conditioning}) in this case.

We will start in Section \ref{sec:LA} with showing how one could \textit{in principle}
form the reduced system using classical linear algebraic techniques. This path connects
our treatment to standard domain decomposition techniques. We then describe an alternative,
more computationally efficient, path in Sections \ref{sec:analytic} and \ref{sec:HPS}.

\begin{figure}[H]
\begin{center}
\begin{subfigure}{.3\linewidth}
    \centering
    \textit{Computational grid}
    \includegraphics[width=\linewidth]{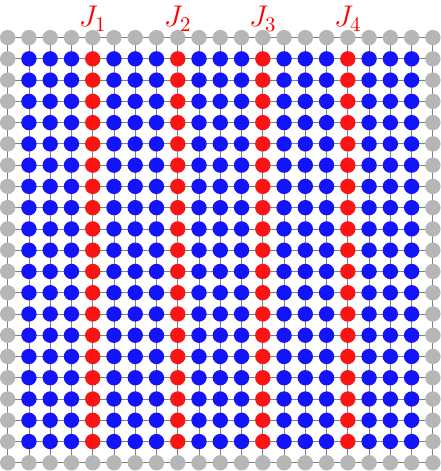}
\end{subfigure}
\hfill
\begin{subfigure}{.3\linewidth}
    \centering
    \textit{Sparsity pattern of }$\mtx{T}$
    \vspace{9pt}
    
    \includegraphics[width=\linewidth]{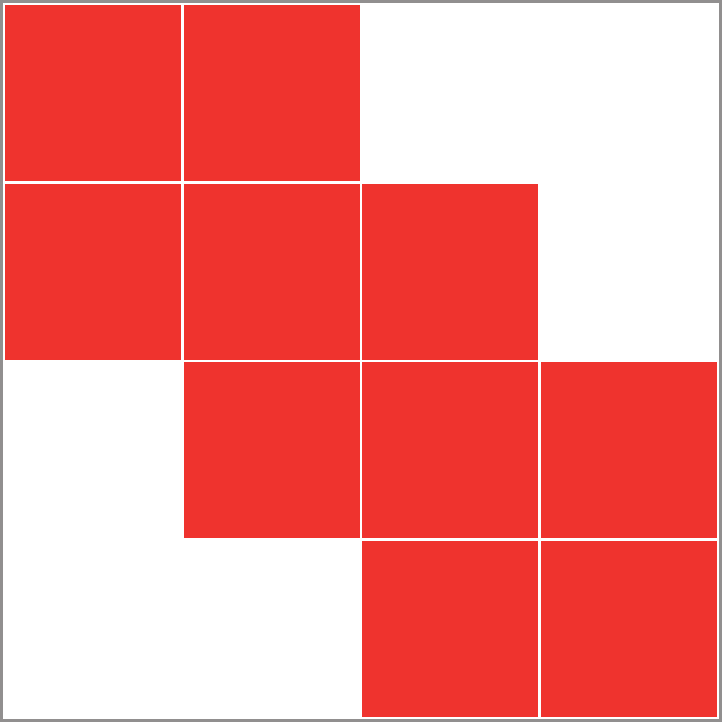}
\end{subfigure}
\hfill
\begin{subfigure}{.3\linewidth}
    \centering
    \textit{Sparsity pattern of }$\mtx{S}$
    \vspace{9pt}
    
    \includegraphics[width=\linewidth]{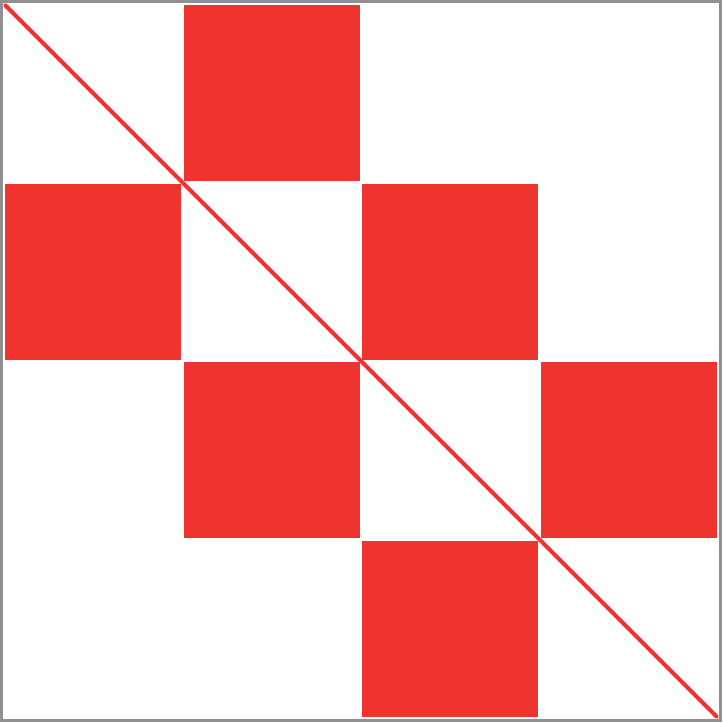}
\end{subfigure}
\end{center}
\caption{The model problem considered in Section \ref{sec:LA}:
A linear system $\mtx{A}$ results from discretization on a computational grid
(gray nodes hold Dirichlet data and are not `active').
The domain is split into thin slabs, separated
by the nodes in $J:=\{J_{j}\}_{j=1}^{4}$ (red). The reduced linear system (middle, see (\ref{eq:modelT})) has the Schur complement coefficient matrix $\mtx{T} = \mtx{A}(J,J) - \mtx{A}(J,J^c)\,\mtx{A}(J^c,J^c)^{-1}\,\mtx{A}(J^c,J)$. The matrix $\mtx{S}$ from (\ref{eq:modelS}) shown on the right is obtained by block diagonal preconditioning of $\mtx{T}$.}
\label{fig:FD}
\end{figure}

\subsection{A numerical derivation}
\label{sec:LA}

For purposes of illustration, let us revisit the toy problem introduced in
Section \ref{sec:toyproblem} where we split the unit square $\Omega = [0,1]^{2}$
into five thin strips, and then discretize (\ref{eq:basic}) using the five-point finite difference
stencil on a regular grid such as the one shown in Figure \ref{fig:FD}.
This results in a linear system
$$
\mtx{A}\vct{u} = \vct{b}
$$
where $\mtx{A}$ is a sparse matrix, where $\vct{u}$ holds the values of the
approximate solution at the interior grid points, and where $\vct{b}$ holds the information
from the boundary condition $f$ and the body load $g$.

To derive the first reduced linear system (\ref{eq:modelT}), we collect the
indices in the mesh separators into the index vector
$$
J = J_{1} \cup J_{2} \cup J_{3} \cup J_{4},
$$
cf.~Figure \ref{fig:FD}.
All remaining interior indices are listed in the complement $J^{\rm c}$.
Performing one step of block Gaussian elimination, we eliminate all nodes
in the interiors of the slabs. The resulting matrix $\mtx{T}$ in (\ref{eq:modelT})
is then simply the Schur complement
$$
\mtx{T} = \mtx{A}(J,J) - \mtx{A}(J,J^c)\,\mtx{A}(J^c,J^c)^{-1}\,\mtx{A}(J^c,J).
$$
The matrix $\mtx{T}$ is block tridiagonal, as any two interfaces $\Gamma_{j}$
and $\Gamma_{j'}$ are disconnected if $|j-j'| > 1$.

As spelled out in Section \ref{sec:toyproblem}, the linear system (\ref{eq:modelS})
is obtained by simply applying block diagonal pre-conditioning to (\ref{eq:modelT}).
In other words, the blocks in $\mtx{S}$ are given by
\begin{equation}
\label{eq:efewin}
\mtx{S}_{j,j'} = -\mtx{T}_{j,j}^{-1}\mtx{T}_{j,j'} = -\mtx{R}_{j}\mtx{A}(I_j,I_j)^{-1}\mtx{A}(I_j,J_{j'})
\end{equation}
in which $I_j$ are the DOFs internal to $\Psi_j$, and $\mtx{R}_j$ is the restriction to $J_j$ (viewed as a subset of $I_j$). From the second part of equation~(\ref{eq:efewin}) we see that $\mtx{S}_{j,j'}$ behaves like a solution operator, followed by a restriction operator. If the blocks of $\mtx{S}$ are to be built explicitly (i.e. we do not use a matrix-free approach), we could in principle do this by first forming $\mtx{T}$, and then
evaluating the formula (\ref{eq:efewin}). In practice, this would be very expensive for large scale problems in 3D, as all blocks in $\mtx{T}$ are dense. It is possible to exploit rank deficiencies in the off-diagonal blocks of $\mtx{T}_{j,j}$ and $\mtx{T}_{j,j'}$ and use, e.g., $\mathcal{H}$-matrix algebra, but a key observation of our work is that $\mtx{S}$ is far more
compressible than $\mtx{T}$, so we will avoid ever forming $\mtx{T}$; Section \ref{sec:analytic} describes how.

\subsection{An analytic derivation}
\label{sec:analytic}

\begin{figure}
\centering
\begin{minipage}{.69\linewidth}
    \begin{subfigure}[t]{.45\linewidth}
\includegraphics[width=\linewidth]{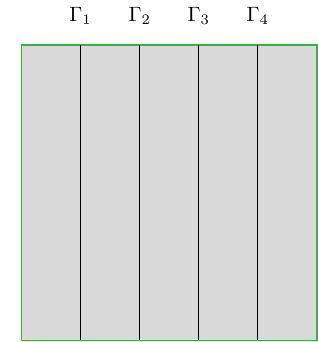}
\caption{}
\end{subfigure}
\begin{subfigure}[t]{.45\linewidth}
    \includegraphics[width=\linewidth]{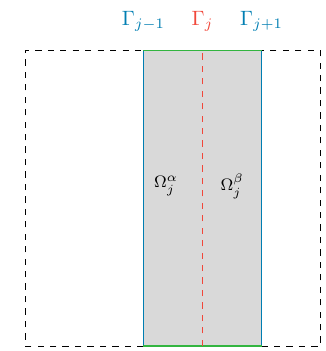}
    \caption{}
\end{subfigure}
\end{minipage}
\begin{minipage}{.3\linewidth}
\textit{(a) BVP on the domain $\Omega$, with known Dirichlet data on $\Gamma$ (green).\\
(b) Local problem on $\Psi_{j}=\Omega_j^{\alpha}\cup\Omega_j^{\beta}$, with known Dirichlet data on $\Gamma\cap\partial\Psi_{j}$. For fixed unknown Dirchlet data $\{u_{j-1},u_{j+1}\}$ on $\Gamma_{j-1}$ and $\Gamma_{j+1}$ (blue), there is a unique
solution $u_{j}$ on $\Gamma_{j}$ (red).}
\end{minipage}
\caption{Continuum domain decomposition described in
Section \ref{sec:analytic}.
}
\label{fig:local_bvp}
\end{figure}

In this section, we take a different path towards deriving the reduced linear system
that starts with a domain decomposition of the continuum problem (\ref{eq:basic}),
\textit{before} discretization. To simplify the discussion, we stick to the simple
toy problem where a square $\Omega = [0,1]^{2}$ has been tessellated into five thin
slabs with interfaces $\{\Gamma_{r}\}_{r=1}^{4}$, as shown
in Figure \ref{fig:local_bvp}(a). We assume that the global problem and each of the local problems is well-posed. Standard techniques for decomposing the full
problem (\ref{eq:basic}) into smaller disconnected subproblems typically involve
enforcing continuity of potentials and normal derivatives across domain boundaries,
and involve forming Dirichlet-to-Neumann operators (or other Poincar\'e-Steklov
operators such as Impedance-to-Impedance maps) for each of the subdomains. A challenge
in this framework is that all boundary operators involve integral or pseudo-integral
operators that have strong singularities at the diagonal. Our objective is to find
an alternative formulation that involves only integral operators with \textit{smooth} kernels.

To simplify the presentation, let us initially consider a problem where the Dirichlet
data on the top and the bottom boundaries are both zero. We will soon return to the
general case.

For the local slab pair, we write
\begin{equation}
\Psi_j = \Omega_j^{(\alpha)} \cup \Omega_j^{(\beta)},
\end{equation}
where $\Omega_j^{(\alpha)}$ and $\Omega_j^{(\beta)}$ are the left and right
halves of the double-wide strip $\Psi_j$.

As a first step, let us for each interface $\Gamma_{j}$ consider a local Dirichlet problem
defined on the double wide strip $\Psi_{j}$ that is enclosed between $\Gamma_{j-1}$ and
$\Gamma_{j+1}$, as shown in Figure \ref{fig:local_bvp}(b).
Suppose that the values of the solution to (\ref{eq:basic}) at the left and the
right boundaries, $u_{j-1}$ and $u_{j+1}$, are known.
Then the solution is uniquely determined everywhere inside $\Psi_{j}$, and in particular
on the line $\Gamma_{j}$.
In other words, there exist linear operators $\mathcal{S}_{j,j-1}$ and $\mathcal{S}_{j,j+1}$ called the \textit{solution operators} such that
$$
u_{j}(x) =
[\mathcal{S}_{j,j-1}u_{j-1}](x) +
[\mathcal{S}_{j,j+1}u_{j+1}](x).
$$
To be precise, if we let $G^{(j)}(x,y)$ denote the Poisson kernel for $\mathcal{A}$ in $\Psi^{(j)}$, then the solution operators take the form
\begin{equation}
\label{eq:debbie1}
[\mathcal{S}_{j,j-1}u_{j-1}](x) = \int_{\Gamma_{j-1}}G^{(j)}(x,y)\,u_{j-1}(y)\,dy,
\end{equation}
and
\begin{equation}
\label{eq:debbie2}
[\mathcal{S}_{j,j+1}u_{j+1}](x) = \int_{\Gamma_{j+1}}G^{(j)}(x,y)\,u_{j+1}(y)\,dy.
\end{equation}
The measure $G^{(j)}(x,y)dy$ is also called the \textit{harmonic measure} of $\mathcal{A}$ in $\Psi^{(j)}$. Provided the Dirichlet Green's functions of $\mathcal{A}$ exist in $\Psi^{(j)}$ (which is always true for symmetric positive definite elliptic PDE operators, see~\cite{DolzmannMuller1995}), the property follows from Green's second identity. If the coefficients of $\mathcal{A}$ are smooth, then, since $x$ and $y$ are never close to each other in (\ref{eq:debbie1}) and (\ref{eq:debbie2}), the Poisson kernels in these integral operators are smooth (see~\cite[Thm.~6.17]{gilbarg2001elliptic}).

\begin{figure}
    \begin{minipage}[t]{.66\linewidth}
    \vspace{0pt}
    \includegraphics[width=.9\linewidth]{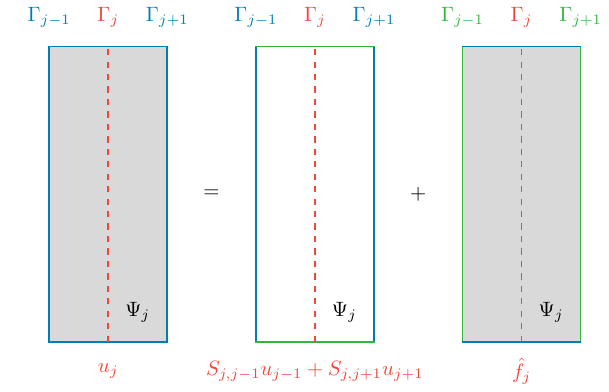}
    \end{minipage}
    \hfill
    \begin{minipage}[t]{.33\linewidth}
    \vspace{5pt}
            \textit{At the central interface $\Gamma_j$ (red), the solution $u_j$ can be written as the sum of the solution operators applied to $\{u_{j-1},u_{j+1}\}$ and the equivalent load $\hat{f}_j$, which is the solution of the local inhomogeneous problem on $\Psi_j$ with zero boundary conditions on $\{\Gamma_{j-1},\Gamma_{j+1}\}$ and the load (gray) of the original global problem.}
            \end{minipage}
    \caption{Illustration of the solution operator principle for general configurations, with nonzero boundary conditions (blue) and load (gray).}
    \label{fig:genBC}
\end{figure}

For the general loaded Poisson problem with nonzero Dirichlet boundary conditions, the solution $u_j$ at interface $\Gamma_j$ can be written as
\begin{equation}\label{eq:equilib1}
    u_j = \mathcal{S}_{j,j-1}u_{j-1}+\mathcal{S}_{j,j+1}u_{j+1}+\hat{f}_j
\end{equation}
where $\hat{f}_j$ is the restriction to $\Gamma_j$ of the solution to the original boundary value problem~(\ref{eq:basic}) restricted to $\Psi_j$, with zero Dirichlet conditions on $\{\Gamma_{j-1},\Gamma_{j+1}\}$. Figure~\ref{fig:genBC} shows a diagram illustrating this principle. Re-writing equation~(\ref{eq:equilib1}) we obtain the \textit{equilibrium equations}
\begin{equation}\label{eq:equilib2}
    -\mathcal{S}_{j,j-1}u_{j-1}+u_j -\mathcal{S}_{j,j+1}u_{j+1}=\hat{f}_j
\end{equation}
which we can accumulate into
\begin{equation}
\label{eq:modelScond}
\mathcal{S}
\left[\begin{array}{rrrr}
u_{1} \\ u_{2} \\ u_{3} \\ u_{4}
\end{array}\right]
=
\left[\begin{array}{rrrr}
\hat{f}_{1} \\ \hat{f}_{2} \\ \hat{f}_{3} \\ \hat{f}_{4}
\end{array}\right],
\qquad\mbox{where}\qquad
\mathcal{S} =
\left[\begin{array}{rrrr}
 \mathcal{I}       & -\mathcal{S}_{1,2} &         &         \\
-\mathcal{S}_{2,1} &  \mathcal{I}       & -\mathcal{S}_{2,3} &        \\
 & -\mathcal{S}_{3,2} &  \mathcal{I}       & -\mathcal{S}_{3,4} \\
 & & -\mathcal{S}_{4,3} &  \mathcal{I}
\end{array}\right].
\end{equation}
For this reason $\mathcal{S}$ is called the \textit{equilibrium operator}.
Of course in practice we do not have access to the Green's kernel, but the discretized solution operators $\mtx{S}_{j,j'}$ can be computed as
\begin{equation}\label{eq:defS}
    \mtx{S}_{j,j'} = -\mtx{R}_j\mtx{A}_j(I,I)^{-1}\mtx{A}_j(I,J)\mtx{R}_{j'}^*
\end{equation}
where $\mtx{A}_j$ is the \textit{local stiffness matrix} for $\Psi_j$, i.e., the discretization of the PDE operator in the double-wide slab $\Psi_j$, $I$ and $J$ are its (local) interior and boundary DOFs, and $\mtx{R}_j$ and $\mtx{R}_{j'}$ denote the restrictions (locally in $\Psi_j$ and $\partial\Psi_j$ respectively) to interface $\Gamma_j$ and $\Gamma_{j'}$ respectively. To be explicit, $\mtx{A}_j$ can be derived from a global discretization, but in our implementation it will be constructed separately for each double-wide slab $\Psi_j$.
\begin{remark}
In principle, even if the original boundary value problem~(\ref{eq:basic}) is not degenerate, it can still happen that one of the local problems is either degenerate or highly unstable. In practice, we have never observed this to happen despite extensive numerical experiments. However, a code targeted at oscillatory problems could be designed to either use other boundary conditions (e.g.~impedance), or to deploy an explicit numerical check for ill-posedness and a mechanism to rectify this as needed.
\end{remark}

\subsection{Discretization of the local problem}
\label{sec:HPS}
The local Dirichlet problem described in Section \ref{sec:analytic} can in
principle be solved with a wide variety of discretization techniques and elliptic
solvers. Our formulation is therefore not tied to the particular local discretizations (and corresponding local stiffness matrices $\{\mtx{A}_j\}_j$), and not tied to the chosen local solvers. All that is required is that the local discretizations are suitably \textit{compatible}. By this we mean that the discretizations of $\Gamma_j$ in each of the double-wide slabs that intersects with $\Gamma_j$ (that is, $\Psi_{j-1}$, $\Psi_{j}$ and $\Psi_{j+1}$) are the same. In principle this requirement can be relaxed through the use of interpolation operators between the various discretizations, but we will not consider this in our manuscript. The $\mathcal{H}$-matrix compression technique considered in this article imposes the additional restriction that all interfaces $\{\Gamma_j\}_{j=1}^{N_{\rm ds}}$ are discretized isomorphically (see Section~\ref{sec:HBS}). This requirement can be lifted through the use of more involved hierarchical matrix formats, but for the applications here it is not restrictive.

\begin{figure}
\centering
\setlength{\unitlength}{1mm}
\begin{picture}(130,70)
\put(00,00){\includegraphics[height=65mm]{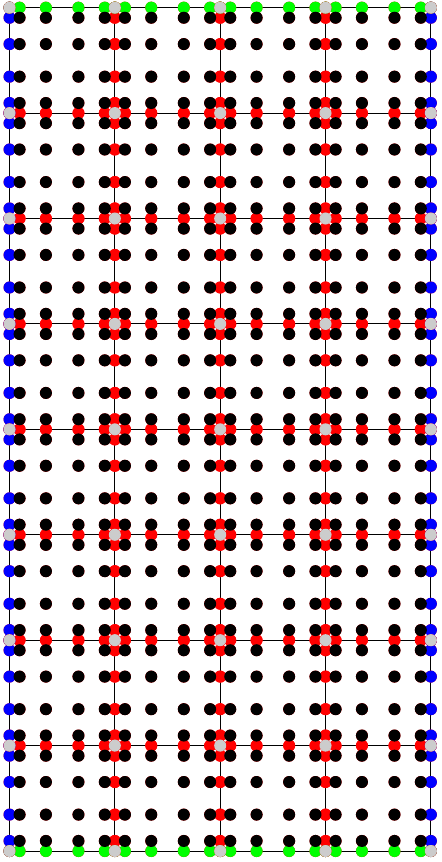}}
\put(38,05){$\Psi_{j}$}
\put(47,34){\begin{minipage}{82mm}
            \textit{The local boundary value problem on each double-wide
            slab $\Psi_{j}$ is solved using a multidomain high
            order spectral method where the slab is tessellated
            into small cells. On each square, a
            Chebyshev grid with $p\times p$ nodes is placed (shown for $p=6$).
            The PDE is enforced directly via spectral differentiation
            and collocation at each node that is interior to a cell (black).
            At each edge node (red), continuity of normal derivatives is
            enforced via spectral differentiation.
            Zero Dirichlet conditions are enforced at the green nodes,
            and general Dirichlet conditions are enforced at the blue nodes.
            Corner nodes (gray) are inactive. See \cite{kump2026hps} for more details.
            }
            \end{minipage}
           }
\end{picture}
\caption{Illustration of the discretization technique described in Section \ref{sec:HPS}
for discretizing the local boundary value problems introduced in Section \ref{sec:analytic}.}
\label{fig:HPS}
\end{figure}

For the high-order experiments in Section~\ref{sec:numerical}, we use a multidomain spectral collocation technique known as \textit{Hierarchical Poincar\'e-Steklov (HPS)}
\cite{2012_spectralcomposite,2013_martinsson_DtN_linearcomplexity} (cf.~also \cite{2003_pfeiffer_spectralmultidomain}). Specifically, we use the implementation from~\cite{kump2026hps}. Following \cite{2022_martinsson_slabLU}, we combine this discretization with a local sparse solver, which is particularly efficient in the present context because the local domains are thin. This HPS implementation is briefly summarized in Figure \ref{fig:HPS}; for further details, see \cite[Ch.~25]{2019_martinsson_fast_direct_solvers}.

\subsection{Rank structure of the off-diagonal blocks}
\label{sec:HBS}

For small two-dimensional problems, one can form the off-diagonal blocks
$\mtx{S}_{j,j-1}$ and $\mtx{S}_{j,j+1}$ in the reduced system (\ref{eq:modelS}) densely,
regardless of the local solver used to approximate the double-slab Dirichlet maps.
However, this dense formation becomes impractical in three dimensions, and even in two
dimensions it can become expensive as the number of interfaces or the local resolution on the subdomains grows.

\begin{figure}[H]
\begin{center}
\setlength{\unitlength}{1mm}
\begin{minipage}{.4\linewidth}
            \textit{Example of a \textit{rank-structured} matrix with weak admissibility.
            Each off-diagonal block (gray) has low numerical rank, and each diagonal block (red) is treated as dense.\\ The tessellation pattern shown is just one example among many possible ones.}


\end{minipage}
\hfill
\begin{minipage}{.45\linewidth}
\includegraphics[height=40mm]{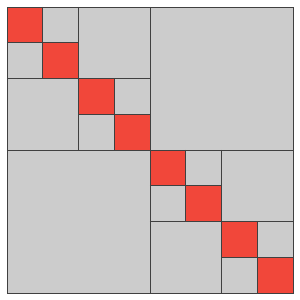}
\end{minipage}
\end{center}
\caption{Illustration of a representative rank-structured matrix, such as an `$\mathcal{H}$-matrix', a `HODLR matrix' as well as an `HBS matrix'.}
\label{fig:tessellation}
\end{figure}

To overcome this problem, we use that the discretized Dirichlet-to-Dirichlet operators inherit exploitable structure from the integral operators (\ref{eq:debbie1}) and (\ref{eq:debbie2}). In particular, the off-diagonal blocks of the $\mtx{S}_{j,j'}$-matrices tend to have very low numerical rank, which means that they can efficiently be represented as \textit{rank-structured matrices}. The idea is to, based on some hierarchical clustering of the DOFs, divide a large dense matrix into a moderate number of blocks in such a way that each block is either of low numerical rank, or is sufficiently small that it can be handled densely. The low-rank matrix blocks correspond to cluster-cluster interactions that are considered `well-separated', in the sense that one expects the local solution kernel from (\ref{eq:debbie1}) and (\ref{eq:debbie2}) to be smooth. Various such hierarchical matrix ($\mathcal{H}$-matrix) formats can be used for this purpose. In this manuscript, we focus on hierarchical matrices with \textit{weak admissibility}. An illustration of a hierarchical matrix with weak admissibility is given in Figure~\ref{fig:tessellation}.
Foundational work in this area was done by Hackbusch and co-workers using the so called
$\mathcal{H}$- and $\mathcal{H}^{2}$-matrix formats
\cite{hackbusch,2008_bebendorf_book,2002_hackbusch_H2,2004_borm_hackbusch,2010_borm_book}.
However, we will use a simpler format, that for our set-up is faster and more efficient:
\textit{Hierarchically Block Separable (HBS)} matrices (sometimes referred to
as \textit{Hierarchically Semi Separable (HSS)})
\cite{2005_martinsson_fastdirect,2012_martinsson_FDS_survey,2009_xia_superfast,2013_xia_randomized,2010_gu_xia_HSS,2016_kressner_rankstructured_review}. 

We will not go into detail here on exactly how HBS matrices are constructed, but a detailed description can be found in \cite{fastHBS}. For the purposes of this text, it suffices to know that an HBS compression of a block $\mtx{S}_{j,j\pm1}\in \mathbb{R}^{N\times N}$ of hierarchical rank $k$ can be computed from $2s = \mathcal{O}(k)$ samples $\mtx{Y}=\mtx{S}_{j,j\pm1}\mtx{\Omega}$ and $\mtx{Z}=\mtx{S}_{j,j\pm1}^*\,\mtx{\Psi}$, where $\mtx{\Omega},\mtx{\Psi}\in \mathbb{R}^{N\times s}$ are Gaussian i.i.d. matrices. Once constructed, the memory footprint of an HBS matrix and the cost for one matrix-vector apply are both $\mathcal{O}(kN)$. The HBS format we use recursively builds a hierarchical approximation based on a \textit{shared cluster tree} for both $\Gamma_j$ and $\Gamma_{j'}$; this is somewhat restrictive, requiring the DOF sets of $\Gamma_{j}$ and $\Gamma_{j'}$ to be isomorphic. This requirement can be lifted by moving to more complicated Hierarchical matrix formats, but for the problems considered here the HBS format is not restrictive.

\begin{figure}[H]
\begin{center}
\setlength{\unitlength}{1mm}
\begin{minipage}{.4\linewidth}
\includegraphics[width=\linewidth]{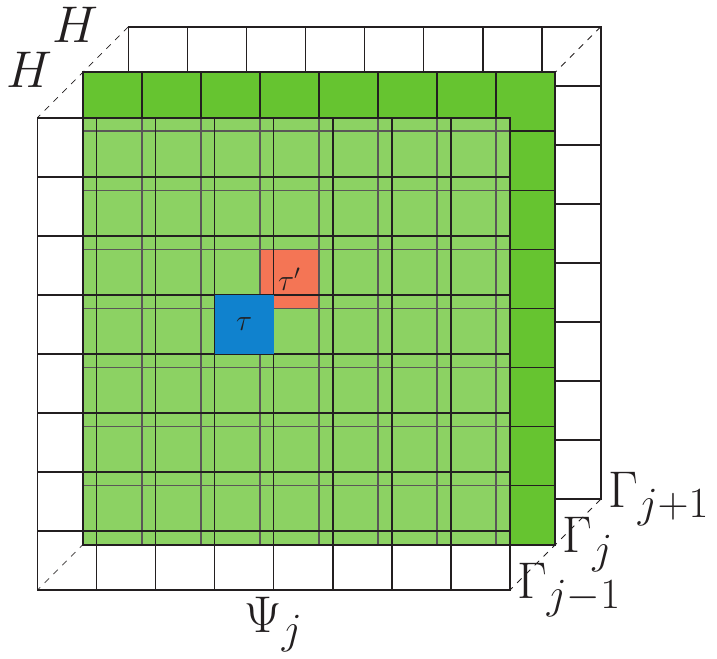}
\end{minipage}
\hspace{1em}
\begin{minipage}{.42\linewidth}
            \textit{Illustration of weak admissibility for the operator $\mtx{S}_{j,j-1}$, defined in the double-wide slab $\Psi_j$ at slab separation $H$. The admissible blocks corresponding to the far-field of the cluster $\tau$ are illustrated in green. For weak admissibility these are the complement of $\tau'$ (red), which is the cluster on interface $\Gamma_j$ that is closest to $\tau$. Note that for the $\mtx{S}$-formulation $\tau$ and $\tau'$ are still separated.}
\end{minipage}

\end{center}
\caption{Illustration of the far field of a cluster $\tau$ with weak admissibility in the $\mtx{S}$-formulation.}
\label{fig:admissibility}
\end{figure}

We stress again that in this manuscript we compress blocks that correspond to \textbf{separated surfaces of source and target points}, $\Gamma_{j'}$ and $\Gamma_j$. This means that our definition of weak admissibility aligns algebraically with the usual definition of weak admissibility for HBS matrices (i.e., only the diagonal blocks of the $\mtx{S}_{j,j'}$ are considered dense), but geometrically it differs slightly from the usual definition (see~\cite[Ch.~15]{2019_martinsson_fast_direct_solvers}). Instead of cluster interactions being considered inadmissible only if the clusters are adjacent, we consider them to be inadmissible if they are of minimal distance, as shown in Figure~\ref{fig:admissibility}. For the 3D case, to our knowledge, this manuscript presents the first use of the randomized HBS compression method from~\cite{fastHBS} for surfaces in 3D.

\begin{figure}[H]
\begin{center}
\setlength{\unitlength}{1mm}
\begin{minipage}[t]{.3\linewidth}
\vspace{0pt}
\includegraphics[width=\linewidth]{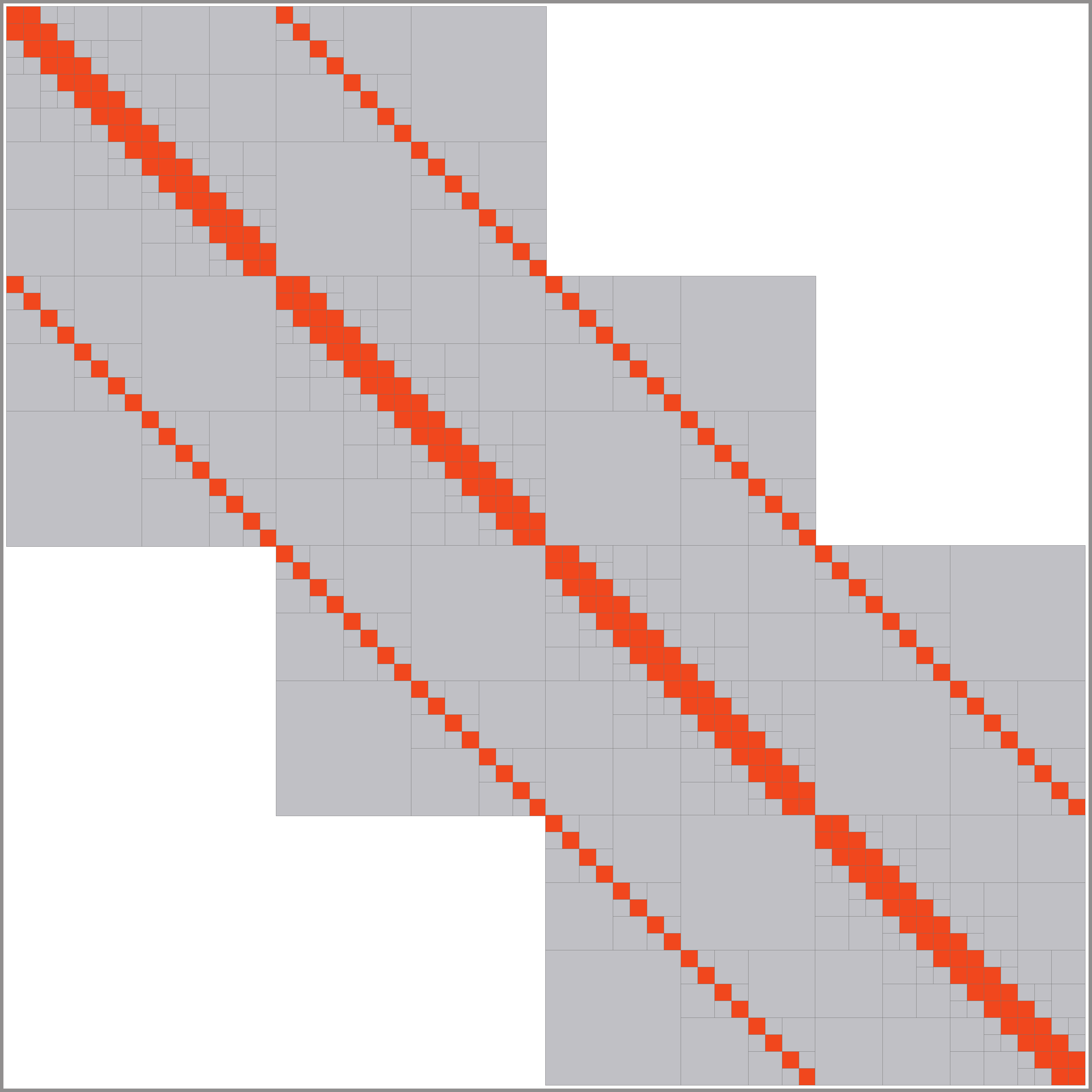}
\end{minipage}
\hfill
\begin{minipage}[t]{.3\linewidth}
\vspace{0pt}
\includegraphics[width=\linewidth]{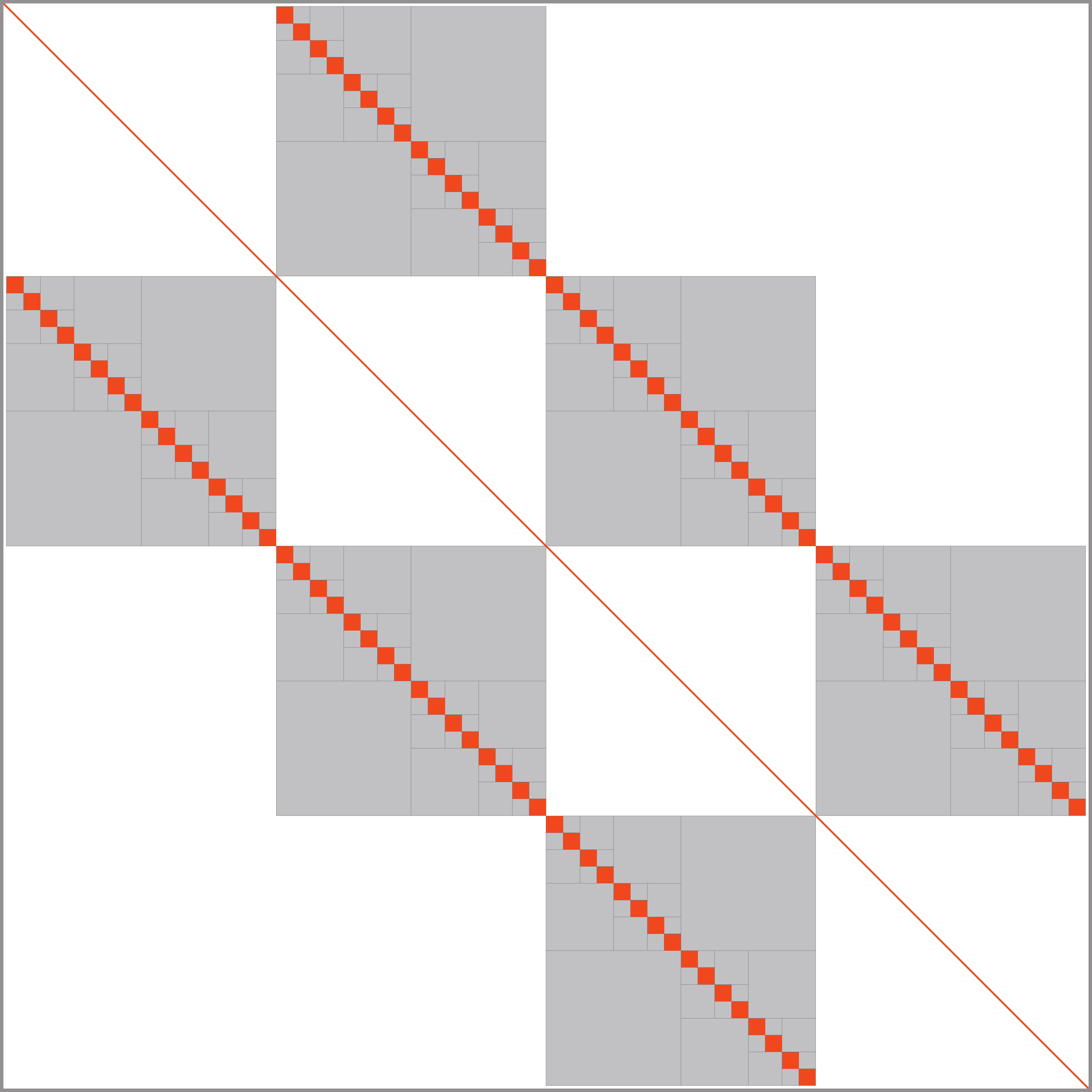}
\end{minipage}
\hfill
\begin{minipage}[t]{.38\linewidth}
\vspace{0pt}
            \textit{The structures of the equilibrium systems $\mtx{T}$ (left) and $\mtx{S}$ (right), using weak admissibility for HBS compression of the off-diagonal blocks. The diagonal blocks of $\mtx{T}$ require strong admissibility. Gray blocks correspond to admissible cluster interactions, stored as low-rank matrices. At fixed rank, the memory footprint for $\mtx{S}$ is much lower than that of $\mtx{T}$.}
\end{minipage}

\end{center}
\caption{The structure of the global equilibrium systems $\mtx{T}$ and $\mtx{S}$ with hierarchical compression.}
\label{fig:globalCompressedS}
\end{figure}

\begin{remark}
The effectiveness of basic HBS compression for $\mtx{S}_{j,j\pm 1}$ crucially relies on the validity of \textit{weak admissibility} (see \cite[Ch.~15]{2019_martinsson_fast_direct_solvers}). This assumption is not always valid, especially in 3D, as the slab width $H$ tends to zero. Indeed, it is known in $\mathcal{H}$-matrix techniques for the Helmholtz equation (see, e.g., \cite{Dirckx}) that for the discretizations of boundary integral operators weak admissibility does not suffice in the high-frequency regime. We return to this in Section~\ref{sec:complexity}.
\end{remark}
With local compression enabled, the global equilibrium systems $\mtx{T}$ and $\mtx{S}$ have the structures shown in Figure~\ref{fig:globalCompressedS}. Here, weak admissibility can be used for the off-diagonal blocks of $\mtx{T}$ and $\mtx{S}$ since these correspond to interactions between separated interfaces, but since the block diagonal of $\mtx{T}$ consists of same-interface DtN maps, weak admissibility no longer suffices there (see also Figure~\ref{fig:gmres_comparison}). For these blocks, strong admissibility is required. Under strong admissibility, a block is treated as low rank only when the corresponding source and target clusters satisfy a prescribed geometric separation condition relative to their diameters. Especially in 3D, this significantly increases the memory footprint and computational complexity required to compress these blocks.
\begin{remark}
As already mentioned in Section~\ref{sec:ConstructS}, one alternative to forming the blocks of $\mtx{T}$ or $\mtx{S}$ explicitly is the ``matrix-free'' approach: they are applied using only the factorized interior solver $\mtx{A}_j$ in $\Psi_j$. This is the standard way of constructing iterative solvers based on substructuring domain decompositions (see, e.g.,~\cite{toselliWidlundDD}). We return to this topic in Section~\ref{sec:complexity}.    
\end{remark}

\subsection{Summary of the Proposed Method}
We now summarize our procedure for constructing the discrete solution maps $\mtx{S}_{j,j'}$ (\textit{local assembly}) and the discretized equilibrium operator $\mtx{S}$ (\textit{global assembly}). The ingredients for our procedure are:
\begin{enumerate}
    \item a domain $\Omega$ in $\mathbb{R}^2$ or $\mathbb{R}^3$,
    \item overlapping double-wide slabs $\{\Psi_j\}_{j=1}^{N_{\text{ds}}}$ covering $\Omega$,
    \item local stiffness matrices $\{\mtx{A}_j\}_{j=1}^{N_{\text{ds}}}$ (these can be computed `on the fly').
\end{enumerate}
With $I_j$ and $J_j$ the local interior and boundary DOFs in $\Psi_j$, set $C_j\subset I_j$ and $J_{j,j'}\subset J_j$ to correspond to the central interface $\Gamma_j$ and the interfaces $\Gamma_{j'}$ on the boundary of $\Psi_j$ respectively.\footnote{As mentioned in Section~\ref{sec:HBS}, for our implementation it is important that the discretization strategy is such that the discretizations of $\Psi_j$ and $\Psi_{j'}$ agree on $\Gamma_j$ and $\Gamma_{j'}$ if $|j-j'|=1$.} The solution map is given by
$$\mtx{S}_{j,j'} = -\mtx{R}_{C_j}\mtx{A}_j(I_j,I_j)^{-1}\mtx{A}_j(I_j,J_{j,j'}).$$
Of course inverting the interior stiffness matrix in this way is computationally undesirable, so instead $\mtx{A}_j(I_j,I_j)$ is factorized into $\mtx{A}_j(I_j,I_j) = \mtx{L}_j\mtx{U}_j$\footnote{We omit possible pivoting here and remark that if the discretization conserves symmetry, we can even factorize $\mtx{A}_j(I_j,I_j)$ using a (pivoted) Cholesky factorization.}. As spelled out in Section~\ref{sec:HBS}, we do not compute $\mtx{S}_{j,j'}$ densely, but approximate it using HBS compression, implemented using the method from~\cite{fastHBS}.

The local and global assembly are summarized in Algorithm~\ref{alg:localS} and~\ref{alg:globalS}.

\begin{algorithm}
\caption{Local $\mtx{S}_{j,j'}$-matrix construction}
\algrenewcommand\algorithmicrequire{\textbf{Input}}
\algrenewcommand\algorithmicensure{\textbf{Output}}
\begin{algorithmic}[1]
\Require interface DOF sets $C_{j}\subset I_j$, $J_{j,j'}\subset J_j$, stiffness matrix $\mtx{A}_j$, factorization $\mtx{A}_j(I_j,I_j)=\mtx{L}_j\mtx{U}_j$, rank and tree $(k_j,\mathcal{T}_j)$
\Ensure Solution operator $\mtx{S}_{j,j'}$
\State $n_j\gets |C_j|$
\Comment{We assume $|J_{j,j'}|=n_j$}
\State $s\gets \alpha k_j+q$
\Comment{sample count; typically $q=10$, $\alpha = 3$}
\State draw $\mtx{\Omega},\mtx{\Psi}\sim \mathcal{N}(0,1)^{n_j\times s}$
\State $\mtx{Y}\gets -\mtx{R}_{C_j}\mtx{A}_j^{-1}\mtx{A}_j(I_j,J_{j,j'})\mtx{\Omega}$\Comment{Computed using the LU-factorization of $\mtx{A}_j$}
\State $\mtx{Z}\gets -\mtx{A}_j(I_j,J_{j,j'})^*\mtx{A}_j^{-*}\mtx{R}_{C_j}^*\mtx{\Psi}$
\State $\mtx{S}_{j,j'}\gets HBS(\mtx{\Omega},\mtx{\Psi},\mtx{Y},\mtx{Z},k_j,\mathcal{T}_j)$
\Comment{HBS compression from~\cite{fastHBS}}
\end{algorithmic}
\label{alg:localS}
\end{algorithm}

\begin{algorithm}
\caption{Total $\mtx{S}$-matrix construction}
\algrenewcommand\algorithmicrequire{\textbf{Input}}
\algrenewcommand\algorithmicensure{\textbf{Output}}
\begin{algorithmic}[1]
\Require\,\,\, Local double slabs $\{\Psi_j\}_{j=1}^{N_{\rm ds}}$, HBS ranks and trees $\{(k_j,\mathcal{T}_j)\}_{j=1}^{N_{\text{ds}}}$
\Ensure Discrete equilibrium operator $\mtx{S}$
\State $\mtx{S} \gets \mtx{I}\in \mathbb{R}^{|J_{\Gamma}| \times |J_{\Gamma}|}$
\For{$j=1,\ldots,N_{\text{ds}}$}
\State compute local discretization $\mtx{A}_j$
\State factorize $\mtx{A}_j(I_j,I_j)=\mtx{L}_j\mtx{U}_j$
\State compute central interface DOFs $C_{j}\subset I_j$
\For{$j'\in \{j-1,j+1\}\cap\{1,\ldots,N_{\text{ds}}\}$}
\State compute interface DOFs $J_{j,j'}\subset J_j$
\State compute HBS approx. $\mtx{S}_{j,j'}$ using Algorithm~\ref{alg:localS}
\State set corresponding block in $\mtx{S}$ to $-\mtx{S}_{j,j'}$
\EndFor
\EndFor
\end{algorithmic}
\label{alg:globalS}
\end{algorithm}

\section{Continuum trace formulation and stability under approximation}
\label{sec:continuum_trace}
The construction in Section~\ref{sec:analytic} is most naturally viewed at the continuum level. The continuous equilibrium system solves for the unknown solution traces on the internal interfaces, and its off-diagonal blocks are solution-transfer operators in the overlapping double-wide slabs.  In this section we first record the fixed-$H$ Fredholm structure of the resulting trace equation.  We then show that standard finite-dimensional trace approximations inherit its stability, and finally account for errors introduced by numerical local solves and rank-structured compression.  Symmetry and positive definiteness are not required in this section, as opposed to in Section~\ref{sec:conditioning}, where the conditioning of the equilibrium system is discussed.

Throughout this section, the slab width $H>0$ is fixed.  We assume that the global boundary value problem and the local Dirichlet problems on the single and double-wide slabs are uniquely solvable.  We also assume sufficient regularity, i.e., that the local Dirichlet solution kernels are sufficiently smooth away from the diagonal.  For uniformly elliptic second-order differential operators with smooth coefficients this is always satisfied. For indefinite problems, it is additionally required that neither the global problem nor any local problem is resonant.

\subsection{Fredholm second-kind structure}
Let $\Gamma_1,\ldots,\Gamma_{N_{\rm ds}}$ denote the internal interfaces, and let $\Psi_j$ be the double slab with central interface $\Gamma_j$.  We use the trace space
  $X=\bigoplus_{j=1}^{N_{\rm ds}} X_j$ with $X_j=L^2(\Gamma_j)$.
The natural Dirichlet trace space is $H^{1/2}(\Gamma_j)$, but the $L^2$ setting is convenient here because separation of source and target interfaces makes the local solution maps smoothing.

As outlined in Section~\ref{sec:analytic}, the local Dirichlet problems define solution-transfer operators
\[
  \mathcal S_{j,j-1}:X_{j-1}\to X_j,
  \qquad
  \mathcal S_{j,j+1}:X_{j+1}\to X_j.
\]
In the present Dirichlet-trace formulation, these are Dirichlet-to-Dirichlet maps between separated interfaces.  Including body loads and physical boundary data gives the row-wise equilibrium equations
\begin{equation}
  -\mathcal S_{j,j-1}u_{j-1}+u_j-\mathcal S_{j,j+1}u_{j+1}
  =\widehat f_j.
  \label{eq:continuum_trace_row}
\end{equation}
The obvious modifications apply at interfaces adjacent to the physical boundary.  In block form,
\begin{equation}
  \mathcal S_H u=(\mathcal I-\mathcal K_H)u=\widehat f,
  \label{eq:continuum_trace_operator}
\end{equation}
where $\mathcal K_H$ contains the solution-transfer operators as its nearest-neighbor off-diagonal blocks.

\begin{proposition}[Fredholm structure of the trace equation]
\label{prop:continuum_fredholm}
Under the assumptions stated above, $\mathcal K_H$ is Hilbert--Schmidt, and hence compact, on $X$.  Consequently, $\mathcal S_H=\mathcal I-\mathcal K_H$ is Fredholm of index zero.  If the homogeneous global boundary value problem is uniquely solvable, then $\mathcal S_H$ is invertible.
\end{proposition}

\begin{proof}
By equations~\eqref{eq:debbie1} and~\eqref{eq:debbie2}, each nonzero block of $\mathcal K_H$ is an integral operator whose source variable lies on $\Gamma_{j\pm1}$ and whose target variable lies on $\Gamma_j$.  These interfaces are separated by a positive distance for fixed $H$.  The corresponding solution kernel is therefore smooth on the compact product $\Gamma_j\times\Gamma_{j\pm1}$ and, in particular, is square integrable.  Each block $\mathcal S_{j,j\pm1}$ is thus Hilbert--Schmidt.  Since $\mathcal K_H$ has only finitely many nonzero blocks, it is Hilbert--Schmidt and hence compact.  The operator $\mathcal S_H=\mathcal I-\mathcal K_H$ is therefore Fredholm of index zero.

It remains to establish injectivity.  Suppose that $u\in X$ satisfies the homogeneous trace equation $\mathcal S_Hu=0$.  Use the components of $u$ as Dirichlet data to solve the homogeneous local problems on the double slabs.  The equilibrium equations ensure that the trace produced on each central interface $\Gamma_j$ is precisely $u_j$.  Two neighboring double-slab solutions therefore solve the same homogeneous equation on their common single-slab overlap and have the same Dirichlet data on the boundary of that overlap.  Uniqueness of the local Dirichlet problem implies that the two solutions agree there.  The local solutions consequently patch together to form a homogeneous solution of the global boundary value problem.  Global uniqueness implies that this solution vanishes, and hence $u=0$.  Thus $\mathcal S_H$ is injective; since it is Fredholm of index zero, it is invertible.
\end{proof}

The essential spectrum of $\mathcal S_H$ is therefore the singleton $\{1\}$.  Any other spectral values are isolated eigenvalues of finite multiplicity, with possible accumulation only at $1$.  Moreover,
\[
  \mathcal S_H^*\mathcal S_H
  =\mathcal I-\mathcal K_H-\mathcal K_H^*+\mathcal K_H^*\mathcal K_H
  =\mathcal I+\text{compact},
\]
so the singular values away from $1$ can likewise accumulate only at $1$.  This is the fixed-$H$ second-kind interpretation of the equilibrium formulation.

\subsection{Finite-dimensional approximation}
Consider a sequence of successively more finely refined finite-dimensional trace spaces $\{X_{j}\}_{j=1}^{\infty}$.
We assume that these are nested, $X_1\subset X_2\subset\cdots\subset X$,
and that their union is dense in $X$. Let $P_n\in\mathcal{B}(X)$ be the orthogonal projection onto $X_n$. Then $P_n\to\mathcal I$ strongly. Compactness upgrades this convergence to operator-norm convergence after composition with $\mathcal K_H$:
\begin{equation}
  \|\mathcal K_H-P_n\mathcal K_HP_n\|\longrightarrow0.
  \label{eq:compact_projection_convergence}
\end{equation}
Indeed, $P_n\to\mathcal I$ uniformly on the relatively compact image under $\mathcal K_H$ of the unit ball, so $\|(\mathcal I-P_n)\mathcal K_H\|\to0$.  Since $P_n=P_n^*$ and $\mathcal K_H^*$ is also compact,
\[
  \|\mathcal K_H(\mathcal I-P_n)\|
  =\|(\mathcal I-P_n)\mathcal K_H^*\|\to0,
\]
and~\eqref{eq:compact_projection_convergence} follows.

Define the projected trace operator on $X_n$ by
  $\mathcal S_{H,n}
  =\mathcal I_{X_n}-P_n\mathcal K_H\big|_{X_n}$.
For comparison on a common space, extend it to $X$ by setting
  $\overline{\mathcal S}_{H,n}
  =\mathcal I-P_n\mathcal K_HP_n$.
The extended operator acts as the identity on $X_n^\perp$ and restricts to $\mathcal S_{H,n}$ on $X_n$.

\begin{proposition}[Stability of projected trace operators]
\label{prop:projected_trace_stability}
Assume that $\mathcal S_H$ is invertible.  Then, for all sufficiently large $n$, the extended operator $\overline{\mathcal S}_{H,n}$ and the projected operator $\mathcal S_{H,n}$ are invertible.  Moreover,
\begin{equation}
  \|\overline{\mathcal S}_{H,n}^{-1}\|
  \le
  \frac{\|\mathcal S_H^{-1}\|}
  {1-\|\mathcal S_H^{-1}\|\,
  \|\overline{\mathcal S}_{H,n}-\mathcal S_H\|},
  \label{eq:projected_inverse_bound}
\end{equation}
and
\begin{equation}
  \|\mathcal S_{H,n}^{-1}\|
  \le \|\overline{\mathcal S}_{H,n}^{-1}\|.
  \label{eq:restricted_projected_inverse_bound}
\end{equation}
In particular, both families of inverse norms are uniformly bounded for sufficiently large $n$.
\end{proposition}

\begin{proof}
By~\eqref{eq:compact_projection_convergence},
\[
  \|\overline{\mathcal S}_{H,n}-\mathcal S_H\|
  =\|\mathcal K_H-P_n\mathcal K_HP_n\|\to0.
\]
For sufficiently large $n$, the product of this norm with $\|\mathcal S_H^{-1}\|$ is smaller than one.  A Neumann-series argument then gives invertibility of $\overline{\mathcal S}_{H,n}$ and the bound~\eqref{eq:projected_inverse_bound}.  Since $\overline{\mathcal S}_{H,n}$ is the identity on $X_n^\perp$ and restricts to $\mathcal S_{H,n}$ on $X_n$, the projected operator is invertible as well, and~\eqref{eq:restricted_projected_inverse_bound} follows by restriction.
\end{proof}

We next apply this operator stability result to the projected equation for a fixed right-hand side.  The following corollary separates the effects of approximating the data and approximating the operator.

\begin{corollary}[Convergence of projected solutions]
\label{cor:projected_solution_convergence}
Let $u\in X$ solve $\mathcal S_Hu=f$, and let $u_n\in X_n$ solve
\[
  \mathcal S_{H,n}u_n=P_nf.
\]
Then, for all sufficiently large $n$,
\begin{equation}
  \|u-u_n\|
  \le
  \|\overline{\mathcal S}_{H,n}^{-1}\|
  \left(
  \|(\mathcal I-P_n)f\|
  +\|\overline{\mathcal S}_{H,n}-\mathcal S_H\|\,\|u\|
  \right).
  \label{eq:projected_solution_error}
\end{equation}
In particular, $u_n\to u$ as $n\to\infty$.
\end{corollary}

\begin{proof}
Viewing $u_n$ as an element of $X$, we have
\[
  \overline{\mathcal S}_{H,n}(u_n-u)
  =P_nf-\overline{\mathcal S}_{H,n}u
  =-(\mathcal I-P_n)f
   -(\overline{\mathcal S}_{H,n}-\mathcal S_H)u.
\]
Applying $\overline{\mathcal S}_{H,n}^{-1}$ gives~\eqref{eq:projected_solution_error}.  Strong convergence of $P_n$ and operator-norm convergence of $\overline{\mathcal S}_{H,n}$ then imply $u_n\to u$.
\end{proof}

The projected construction above provides a stable reference approximation of the continuum equation.  A collocation or finite-difference implementation need not be literally obtained by orthogonal projection; instead, the computed operator can be compared with the stable projected operator through a standard perturbation argument (see, e.g.,~\cite[Ch.~3-4]{atkinson1997}).

\begin{remark}[Inexact local solves and compression]
\label{rem:inexact_local_solves}
Proposition \ref{prop:projected_trace_stability} provides a stability result that can be
used to quantify how errors introduced by numerical local solves and rank-structured compression affect the computed solution. To illustrate, let
\[
  \widehat{\mathcal S}_{H,n}
  =\mathcal I_{X_n}-\widehat{\mathcal K}_{H,n}
\]
denote the operator obtained after trace discretization, numerical local Dirichlet solves, quadrature or collocation, and rank-structured compression, and set
\[
  \eta_n
  :=\|\widehat{\mathcal S}_{H,n}-\mathcal S_{H,n}\|.
\]
If
$
  \eta_n\,\|\mathcal S_{H,n}^{-1}\|<1,
$
then the standard Neumann-series perturbation argument shows that $\widehat{\mathcal S}_{H,n}$ is invertible and
\begin{equation}
  \|\widehat{\mathcal S}_{H,n}^{-1}\|
  \le
  \frac{\|\mathcal S_{H,n}^{-1}\|}
  {1-\eta_n\|\mathcal S_{H,n}^{-1}\|}.
  \label{eq:computed_inverse_bound}
\end{equation}
Moreover, if
$\mathcal S_{H,n}u_n=P_nf$
and
$\widehat{\mathcal S}_{H,n}\widehat u_n=\widehat f_n$,
then
\begin{equation}
  \|\widehat u_n-u_n\|
  \le
  \frac{\|\mathcal S_{H,n}^{-1}\|}
  {1-\eta_n\|\mathcal S_{H,n}^{-1}\|}
  \left(
    \|\widehat f_n-P_nf\|+\eta_n\|u_n\|
  \right).
  \label{eq:computed_solution_error}
\end{equation}
The point here is not to estimate the individual error terms; rather, (\ref{eq:computed_solution_error}) shows how they control the difference between the computed solution and the stable projected solution, once the local solves and compression attain sufficient accuracy.
\end{remark}

The results in this section concern refinement of a fixed continuum trace equation at fixed physical slab width $H$.  They do not quantify how $\|\mathcal S_H^{-1}\|$ depends on $H$, on the frequency in an oscillatory problem, or on proximity to local or global resonances.  
Section~\ref{sec:conditioning} addresses the $H$-dependence of the conditioning for compatible symmetric positive-definite discretizations.
In that setting, the locally constructed solution-transfer equation is algebraically equivalent to a block-scaled Schur complement and admits a quantitative energy estimate.

\section{Compatible SPD discretizations}
\label{sec:conditioning}
Sections~\ref{sec:LA} and~\ref{sec:analytic} introduced two routes to the reduced interface equation.  Conceptually, we regard the continuum-first route of Section~\ref{sec:analytic} as primary: first partition the domain into overlapping double slabs and define the local solution-transfer operators, and then approximate those operators using suitable local discretizations. This is in contrast to the ``discretize first, then eliminate degrees algebraically by forming Schur complements'' approach of Section \ref{sec:LA}. In many standard settings, the order of these operations can be reversed without changing the resulting discrete interface equation.  This occurs, for example, for standard low-order finite differences on matching grids and for conforming finite-element discretizations with matching trace spaces, when the local discrete spaces and operators are obtained by restriction from a single global discretization.  One may then first discretize the global problem, eliminate the slab interiors, and form a block-scaled Schur complement.  We refer to such a realization as a \textit{compatible} discretization.

Throughout this section, we restrict attention to the case in which the boundary value problem (\ref{eq:basic}) is symmetric and coercive, and the compatible discretization matrix is symmetric positive definite (SPD).  The continuum analysis in Section~\ref{sec:continuum_trace} does not require this restriction or the equivalence described above.  The purpose of the present section is to exploit the compatible SPD setting, where it supplies a natural discrete energy structure and a quantitative estimate of the dependence on the slab width $H$.

\subsection{A compatible Schur-complement realization}
Let $\mtx{A}\vct{u}=\vct{b}$ be the linear system obtained from a compatible SPD discretization of the boundary value problem (\ref{eq:basic}).  Partition the active degrees of freedom into the interface set $J$, corresponding to $\Gamma_1,\ldots,\Gamma_{N_{\rm ds}}$, and the remaining slab-interior set $I$.
(In Figure \ref{fig:FD}, $J$ are the red nodes, and $I$ are the blue nodes.)
Boundary values prescribed on the physical boundary are understood to have been moved to the right-hand side.  Eliminating the interior degrees of freedom gives the interface Schur complement
\begin{equation}
  \mtx{T}_H
  =\mtx{A}_{JJ}-\mtx{A}_{JI}\mtx{A}_{II}^{-1}\mtx{A}_{IJ}.
  \label{eq:T_def_theory}
\end{equation}
Since $\mtx{A}$ is SPD, so is $\mtx{T}_H$,
and so is each of its principal diagonal blocks.  Split $J=J_1\cup\cdots\cup J_{N_{\rm ds}}$ by interface, and define
\begin{equation}
  \mtx{D}_H
  =\operatorname{diag}\bigl(
    \mtx{T}_{11},\ldots,
    \mtx{T}_{N_{\rm ds},N_{\rm ds}}
  \bigr).
  \label{eq:D_def_theory}
\end{equation}
The corresponding block-scaled interface matrix is
\begin{equation}
  \mtx{S}_H=\mtx{D}_H^{-1}\mtx{T}_H.
  \label{eq:S_def_theory}
\end{equation}
For a slab decomposition, $\mtx{T}_H$ is block tridiagonal.  Compatibility of the global and local discretizations implies that the off-diagonal blocks
$
  -\mtx{T}_{jj}^{-1}\mtx{T}_{j,j\pm1}
$
are exactly the discrete solution-transfer operators obtained by solving the corresponding local Dirichlet problems on the double slab $\Psi_j$.  Thus~\eqref{eq:S_def_theory} is a discrete realization of the continuum trace equation from Section~\ref{sec:continuum_trace}; the Schur-complement construction is a tool for analysis, rather than a prerequisite for computing the method.

The Schur complement has the variational characterization (see, e.g.,~\cite{SmithSform})
\begin{equation}
  \vct{v}^T\mtx{T}_H\vct{v}
  =\min_{\vct{w}(J)=\vct{v}}
  \vct{w}^T\mtx{A}\vct{w}.
  \label{eq:schur_energy_min}
\end{equation}
Hence $\vct{v}^T\mtx{T}_H\vct{v}$ is the energy of the global discrete $\mtx{A}$-harmonic extension of the interface data $\vct{v}$.  Similarly,
\[
  \vct{v}^T\mtx{D}_H\vct{v}
  =\sum_{j=1}^{N_{\rm ds}}
    \vct{v}_j^T\mtx{T}_{jj}\vct{v}_j
\]
is the corresponding additive local interface energy.  The conditioning estimate below follows from comparing these two energies.

\subsection{Energy scaling and red-black structure}
Although $\mtx{S}_H=\mtx{D}_H^{-1}\mtx{T}_H$ is generally not symmetric in the Euclidean inner product, it is symmetric positive definite in the natural $\mtx{D}_H$-inner product.

\begin{proposition}[Energy symmetry and symmetric scaling]
\label{prop:energy_symmetry}
Let $\mtx{T}_H$ and $\mtx{D}_H$ be defined by~\eqref{eq:T_def_theory} and~\eqref{eq:D_def_theory}.  Then $\mtx{S}_H=\mtx{D}_H^{-1}\mtx{T}_H$ is self-adjoint and positive definite in the inner product
\[
  \langle \vct{x},\vct{y}\rangle_{\mtx{D}_H}
  =\vct{x}^T\mtx{D}_H\vct{y}.
\]
Equivalently, $\mtx{S}_H$ is similar to the Euclidean SPD matrix
\begin{equation}
  \widetilde{\mtx{S}}_H
  =\mtx{D}_H^{-1/2}\mtx{T}_H\mtx{D}_H^{-1/2}.
  \label{eq:S_tilde}
\end{equation}
\end{proposition}

\begin{proof}
For any vectors $\vct{x}$ and $\vct{y}$,
$
  \langle \mtx{S}_H\vct{x},\vct{y}\rangle_{\mtx{D}_H}
  =\vct{x}^T\mtx{T}_H\vct{y}
  =\langle \vct{x},\mtx{S}_H\vct{y}\rangle_{\mtx{D}_H},
$
and
$
  \langle \mtx{S}_H\vct{x},\vct{x}\rangle_{\mtx{D}_H}
  =\vct{x}^T\mtx{T}_H\vct{x}>0
$
for every nonzero $\vct{x}$.  The similarity relation
$
  \mtx{S}_H
  =\mtx{D}_H^{-1/2}\widetilde{\mtx{S}}_H\mtx{D}_H^{1/2}
$
then gives the equivalent Euclidean SPD representation.
\end{proof}

Writing $\|\vct{x}\|_{\mtx{D}_H}^2=\vct{x}^T\mtx{D}_H\vct{x}$ and using the induced operator norm, the natural $\mtx{D}_H$-condition number is
\begin{equation}
  \kappa_{\mtx{D}_H}(\mtx{S}_H)
  :=\|\mtx{S}_H\|_{\mtx{D}_H}
    \|\mtx{S}_H^{-1}\|_{\mtx{D}_H}
  =\kappa_2\bigl(\widetilde{\mtx{S}}_H\bigr).
  \label{eq:energy_condition_number}
\end{equation}
This energy condition number is the quantity controlled below.  

The upper spectral bound has a useful red-black interpretation.  For a nonperiodic slab chain, color the odd-numbered interfaces red and the even-numbered interfaces black.  The same two-coloring applies to a periodic decomposition when $N_{\rm ds}$ is even.  After this ordering,
\[
  \mtx{T}_H
  =\begin{bmatrix}
      \mtx{T}_{\rm rr}&\mtx{T}_{\rm rb}\\
      \mtx{T}_{\rm br}&\mtx{T}_{\rm bb}
    \end{bmatrix},
  \qquad
  \mtx{D}_H
  =\begin{bmatrix}
      \mtx{T}_{\rm rr}&\mtx{0}\\
      \mtx{0}&\mtx{T}_{\rm bb}
    \end{bmatrix},
\]
where $\mtx{T}_{\rm rr}$ and $\mtx{T}_{\rm bb}$ are block diagonal.  Consequently,
\begin{equation}
  \mtx{S}_H
  =\underbrace{
    \begin{bmatrix}
      \mtx{T}_{\rm rr}^{-1}&\mtx{0}\\
      \mtx{0}&\mtx{0}
    \end{bmatrix}\mtx{T}_H}_{\mtx{P}_{\rm r}}
  +\underbrace{
    \begin{bmatrix}
      \mtx{0}&\mtx{0}\\
      \mtx{0}&\mtx{T}_{\rm bb}^{-1}
    \end{bmatrix}\mtx{T}_H}_{\mtx{P}_{\rm b}}.
  \label{eq:red_black_projectors}
\end{equation}
Direct calculation shows that $\mtx{P}_{\rm r}$ and $\mtx{P}_{\rm b}$ are idempotent and self-adjoint in the $\mtx{T}_H$-inner product.  Thus each is a $\mtx{T}_H$-orthogonal projection.  They are not, in general, mutually orthogonal.

\begin{proposition}[Red-black upper bound]
\label{prop:red_black_upper_bound}
For the two-colorable slab decompositions described above,
\[
  \lambda_{\max}\bigl(\widetilde{\mtx{S}}_H\bigr)\le 2.
\]
\end{proposition}

\begin{proof}
Since $\mtx{P}_{\rm r}$ and $\mtx{P}_{\rm b}$ are orthogonal projections in the $\mtx{T}_H$-inner product,
\[
  \|\mtx{P}_{\rm r}\|_{\mtx{T}_H}\le1,
  \qquad
  \|\mtx{P}_{\rm b}\|_{\mtx{T}_H}\le1.
\]
Equation~\eqref{eq:red_black_projectors} therefore gives
\[
  \rho(\mtx{S}_H)
  \le\|\mtx{S}_H\|_{\mtx{T}_H}
  \le2.
\]
The matrices $\mtx{S}_H$ and $\widetilde{\mtx{S}}_H$ have the same positive eigenvalues, so the stated bound follows.
\end{proof}

\subsection{Stable splitting and conditioning}
The lower spectral bound requires the following standard stable-splitting estimate.

\begin{assumption}[Stable splitting]
\label{ass:stable_decomposition}
There exists a constant $C_{\rm ss}$, independent of the local mesh size and the number of degrees of freedom used in each slab, such that
\begin{equation}
  \vct{v}^T\mtx{D}_H\vct{v}
  \le C_{\rm ss}H^{-2}\,
       \vct{v}^T\mtx{T}_H\vct{v}
  \label{eq:stable_decomp_assumption}
\end{equation}
for every interface vector $\vct{v}$.  The constant may depend on the ellipticity constants, the global domain, and the regularity of the slab decomposition.
\end{assumption}

Assumption~\ref{ass:stable_decomposition} is a mild and standard requirement for the compatible discretizations considered here.  It holds, in particular, for standard low-order finite-difference discretizations on regular grids and for conforming low-order finite-element discretizations on shape-regular meshes; see, e.g.,~\cite{smithDD,toselliWidlundDD}.  At the continuum level, estimates of this type are obtained by multiplying the global energy-minimizing extension by a partition of unity subordinate to the overlapping slabs.  The gradients of the cutoff functions scale like $H^{-1}$, producing the factor $H^{-2}$ in the energy estimate.  The corresponding discrete estimates follow from the usual stable-decomposition arguments in additive Schwarz theory.  For high-order Galerkin or spectral-element discretizations, uniformity in the polynomial order additionally requires the relevant interpolation and decomposition estimates to be $p$-stable.

\begin{theorem}[Energy and spectral bounds]
\label{thm:condDisc}
Assume that $\mtx{A}$ is SPD, that the interface graph is two-colorable as above, and that Assumption~\ref{ass:stable_decomposition} holds.  Then
\begin{equation}
  C_{\rm ss}^{-1}H^2
  \le
  \lambda_{\min}\bigl(\widetilde{\mtx{S}}_H\bigr)
  \le
  \lambda_{\max}\bigl(\widetilde{\mtx{S}}_H\bigr)
  \le 2.
  \label{eq:scaled_eigenvalue_bounds}
\end{equation}
Consequently,
\begin{equation}
  \kappa_{\mtx{D}_H}(\mtx{S}_H)
  =\kappa_2\bigl(\widetilde{\mtx{S}}_H\bigr)
  \le 2C_{\rm ss}H^{-2}.
  \label{eq:condition_bound}
\end{equation}
In particular, for fixed physical slab width $H$, this bound is independent of the local discretization resolution.
\end{theorem}

\begin{proof}
For any nonzero vector $\vct{z}$, set $\vct{v}=\mtx{D}_H^{-1/2}\vct{z}$.  Then
\[
  \frac{\vct{z}^T\widetilde{\mtx{S}}_H\vct{z}}
       {\vct{z}^T\vct{z}}
  =
  \frac{\vct{v}^T\mtx{T}_H\vct{v}}
       {\vct{v}^T\mtx{D}_H\vct{v}}.
\]
The stable-splitting estimate~\eqref{eq:stable_decomp_assumption} gives the lower bound
\[
  \frac{\vct{v}^T\mtx{T}_H\vct{v}}
       {\vct{v}^T\mtx{D}_H\vct{v}}
  \ge C_{\rm ss}^{-1}H^2.
\]
The upper bound follows from Proposition~\ref{prop:red_black_upper_bound}.  The condition-number estimate is then immediate from~\eqref{eq:energy_condition_number}.
\end{proof}

Because $\mtx{S}_H$ and $\widetilde{\mtx{S}}_H$ are similar, they have the same eigenvalues. Hence the spectral ratio
\[
\ensuremath{\kappa_{\rho}(\mtx{S}_H):=\frac{\lambda_{\max}(\mtx{S}_H)}{\lambda_{\min}(\mtx{S}_H)}=\kappa_2(\widetilde{\mtx{S}}_H)\le 2C_{\rm ss}H^{-2}}
\]
is also uniform under local refinement. This is a spectral statement and, for a nonnormal matrix, need not equal the Euclidean condition number $\kappa_2(\mtx{S}_H)$. The supplementary material reports weighted spectral diagnostics for representative stencil and HPS collocation realizations. Over the tested range, the ratio $\kappa_2(\mtx{S}_H)/\kappa_{\rho}(\mtx{S}_H)$ remains close to one, suggesting that the spectral bound is indicative of the Euclidean conditioning in these examples. These observations are empirical and are not used in the proof of Theorem~\ref{thm:condDisc}.

\subsection{Collocation and oscillatory problems}
Rather than a compatible SPD discretization, the HPS implementation used in many of the numerical experiments is a strong-form collocation method. Its raw matrices are generally nonsymmetric and nonnormal in the Euclidean inner product. We therefore interpret HPS through the continuum approximation framework of Section~\ref{sec:continuum_trace}; the energy-norm theorem above does not apply directly. The supplementary material reports weighted spectral diagnostics for representative HPS and stencil realizations. In the examples studied, the Euclidean condition number tracks the spectral ratio closely.

For oscillatory problems (where positivity is no longer satisfied), the trace formulation separates two effects.  At fixed frequency and fixed slab width, the interface equation is built from local solution operators rather than from a discretized differential operator, so local refinement does not by itself introduce the large singular values associated with discretizing an unbounded operator. Frequency-dependent stability, including proximity to global or local resonances, remains a separate issue. The Helmholtz experiments investigate the corresponding behavior empirically.


\section{Computational complexity}\label{sec:complexity}
In this section, we analyze the complexity of our proposed solver, and compare it to the ``matrix-free'' and non-overlapping approach.

In what follows, we restrict our attention to the 3D version of our solver, but the same analysis can be applied to the 2D case. The total computational complexity of our proposed method is dominated by the cost for the construction and factorization of the local stiffness matrices on the double-wide slabs $\{\Psi_j\}_{j=1}^{}N_{\rm ds}$. We will denote the selected HBS rank by $k$.

We analyze the case of sparse direct solvers for the interior slab solves, assuming that nested dissection is used. We analyze additional costs associated with the use of high order discretizations in Remark \ref{rem:hps_complexity}.

Consider a 3D domain, discretized with $N$ total DOFs. We assume for simplicity that $N = n_1n_2n_3$, with $n_1$ the number of DOFs in the slab direction, and $n_2=n_3$ the number of DOFs in the $y$ and $z$-axis. There are three primary costs to analyze:
\begin{enumerate}
    \item the cost of factorizing each slab volume,
    \item the cost of compressing $\mtx{S}$ using rank structure, and 
    \item the cost of applying $\mtx{S}$ within GMRES iterations to compute the interface solution. 
\end{enumerate}

\noindent \textbf{Sparse Factorization of slab volumes.}
The domain is divided into $N_{\text{ds}}$ subdomains, and we impose that the slabs are \textit{thin}, such that
$\frac{2n_1}{N_{\text{ds}}} \;\leq\; n_2= n_3$.
Each double-wide slab volume therefore contains $\tfrac{2n_1}{N_{\text{ds}}} \times n_2 \times n_3$ points, leading to costs
$$\text{Factorization:}\ \mathcal O\left( \left(\frac {n_1}{N_{\text{ds}}}\right)^3 n_2^3 \right), \qquad \text{Storage:}\ \mathcal O\left( \left(\frac {n_1}{N_{\text{ds}}}\right)^2 n_2^2 \right), \qquad \text{for each slab}.$$
These complexities are due to the cost of factorizing and storing (respectively) the largest nested dissection separator (see~\cite[Ch.~20]{2019_martinsson_fast_direct_solvers}).
The total cost of factorizing and storing all the volumes' degrees of freedom are
$$\text{Factorization:}\ \mathcal O\left( \frac {n_1^{3/2}}{N_{\text{ds}}^2} N^{3/2}\right), \qquad \text{Storage:}\ O\left( \frac {n_1}{N_{\text{ds}}} N\right), \qquad \text{for all slabs}.$$
When $n_1 = n_2 = n_3$, the complexity costs are
$\mathcal{O}\!\left(\tfrac{N^2}{N_{\text{ds}}^2}\right)$ and
$\mathcal{O}\!\left(\tfrac{N^{4/3}}{N_{\text{ds}}}\right)$
for nested dissection, respectively.

\vspace{0.5em}

\noindent \textbf{Randomized rank-structured compression of $\mtx S$.}
The matrix $\mtx{S}$ is block tridiagonal and acts only on the interfaces.
There are $N_{\text{ds}}$ interfaces, each of size $n_2 n_3=N/ n_1$.
Assume the HBS rank of the off-diagonal blocks is $k$.
Acquiring randomized samples requires $\mathcal{O}(k)$ applications of $\mtx{S}$,
equivalently $\mathcal{O}(k)$ solves with the factorized local stiffness matrices
(see Algorithm~\ref{alg:localS}).
Post-processing these samples is linear in the interface size, and hence sublinear in $N$.
The overall costs are
$$\text{Sampling cost}: \mathcal O\left( k\ \frac {n_1}{N_{\text{ds}}} N\right), \qquad \text{HBS construction}: \mathcal O\left(k^2\ \frac {N_{\text{ds}}}{n_1}\  N\right).$$
For details on the construction, see~\cite{fastHBS}. When $n_1 = n_2 = n_3$, these costs scale as $\mathcal O\left(k\ \frac{N^{4/3}}{N_{\rm ds}}\right)$ and $\mathcal O\left(k^2\ N_{\rm ds}\ N^{2/3}\right)$, respectively.

\vspace{0.5em}

\noindent \textbf{Iteratively solving $\mtx{S}\vct{u}_{\Gamma} = \vct{f}_{\Gamma}$.}
Let $m$ denote the number of Krylov iterations required to reach the prescribed tolerance. For the compatible symmetrically scaled SPD system, Theorem~\ref{thm:condDisc} gives $\kappa_2(\widetilde{\mtx{S}}_H)=\mathcal O(H^{-2})$; when $H\asymp N_{\rm ds}^{-1}$, the standard CG estimate therefore gives the upper bound $m=\mathcal O(N_{\rm ds})$. For HPS-discretized and Helmholtz systems, we leave $m$ as iteration count, but let us mention that for SPD systems discretized using HPS, we observe that $m=\mathcal{O}(H^{-1})$ as well, and for Helmholtz-type systems, the iteration count is also independent of the mesh fineness of the local subdomains.


Once the compressed representation of $\mtx S$ is constructed, one matrix--vector product costs $\mathcal O(kN_{\rm ds}n_2^2)$, while orthogonalization in unrestarted GMRES costs $\mathcal O(m^2N_{\rm ds}n_2^2)$ over the full solve. Consequently,
\[
\text{total solve cost}:\quad
\mathcal O\!\left((k mN_{\rm ds}+m^2N_{\rm ds})n_2^2\right)
=
\mathcal O\!\left(\left(k\frac{mN_{\rm ds}}{n_1}+\frac{m^2N_{\rm ds}}{n_1}\right)N\right).
\]

In the case where $n_2^2$ or $N_{\rm ds}$ are prohibitively large, BiCGStab can be used instead of GMRES to avoid the Krylov orthogonalization cost, at the cost of possibly increasing the iteration count.
\begin{remark}[Weak scaling in the parallel setting]
Algorithm~\ref{alg:globalS} shows that the construction of the global equilibrium operator is embarrassingly parallel.
Once built, matrix–vector products with $\mtx{S}$ can also be parallelized.
If the number of processors and interfaces grow proportionally,
then the cost of sparse factorization, randomized compression,
and application of $\mtx{S}$ all reduce by a factor of $N_{\text{ds}}$.
\label{rem:parallel}
\end{remark}

\begin{remark}[Complexity of the Hierarchical Poincaré–Steklov discretization]
For a fixed HPS tiling, the degrees of freedom grow as
$N = \mathcal{O}(p^3)$.
The general complexity analysis above applies to the HPS discretization as well;
however, one must also account for the additional cost of factorizing the local
differential operators on each subdomain.
These additional costs scale as $\mathcal{O}(p^6 N)$ overall,
as described in detail in \cite{kump2026hps}.
\label{rem:hps_complexity}
\end{remark}
\subsection{Comparison to the matrix-free approach}\label{subsec:compare_to_matfree}
As outlined in Section~\ref{sec:construct_reduced_system}, our compressed overlapping formulation is not the only possible candidate for a slab-based domain decomposition. Two alternatives are:
\begin{enumerate}
    \item per double-wide slab, compute the LU-factors of $\mtx{A}_j$ and apply $\mtx{S}_{j,j\pm 1}$ using local solves, or
    \item form the compressed $\mtx{T}$-system based on the same (single) slabs, using a non-overlapping domain decomposition (imposing continuity of Neumann traces).
\end{enumerate}
These are referred to, respectively, as the \textit{matrix-free} and \textit{non-overlapping} approach. In this section we compare, in terms of asymptotics, the difference in approach between overlapping matrix-free and using HBS compression for the $\mtx{S}$-system. In Section~\ref{subsec:gmres_iters} we compare the overlapping and non-overlapping approach.

In Table~\ref{tab:S_implementation_options} we summarize the contrast between the HBS approach and the matrix-free approach. For simplicity, we restrict our attention to the 3D case with global DOF counts $n_1=n_2=n_3=n$ and set $N = n^3$. We instead write the solve costs in terms of the actual Krylov iteration count $m$, but let us mention that, in all our numerical observations, $m=\mathcal{O}(N_{\rm ds})$
\begin{table}[H]
\centering
\small

\setlength{\tabcolsep}{4pt}
\renewcommand{\arraystretch}{1.2}
\begin{tabular*}{\textwidth}{@{\extracolsep{\fill}}p{0.16\textwidth}p{0.4\textwidth}p{0.4\textwidth}@{}}
\hline
 & Matrix-free $\mtx{S}$ & HBS $\mtx{S}$ \\
\hline
Setup FLOPs & $\mathcal O\left(\frac {N^{2}}{N_{\text{ds}}^2} \right)$\newline $\rightarrow$ $N_{\rm ds}$ local LU factorizations & $\mathcal O\left( \frac {N^2}{N_{\text{ds}}^2} \right) + \mathcal{O}\left(N^{4/3}(k/N_{\rm ds}+k^2N_{\rm ds})\right)$\newline$\rightarrow$ $N_{\rm ds}$ local LU factorizations $+$ $\mathcal{O}(kN_{\rm ds})$ local solves + $N_{\rm ds}$ HBS compressions\\
Storage requirements & $\mathcal{O}\left(\frac{N^{4/3}}{N_{\rm ds}}\right)$ $\rightarrow$ sparse factors. & $\mathcal{O}(kN_{\rm ds}N^{2/3})$ $\rightarrow$ HBS matrices\\
Matrix-vector time & $\mathcal{O}(\frac{N^{4/3}}{N_{\rm ds}})$ $\rightarrow$ $N_{\rm ds}$ local LU solves & $\mathcal{O}(kN_{\rm ds} N^{2/3})$ $\rightarrow$ $N_{\rm ds}$ HBS matvecs\\
Solve FLOPs & $\mathcal O\!\left(mN^{4/3}/N_{\rm ds}+m^2N_{\rm ds}N^{2/3}\right)$ & $\mathcal O\!\left(kmN_{\rm ds}N^{2/3}+m^2N_{\rm ds}N^{2/3}\right)$\\
\hline
\end{tabular*}
\caption{Comparison of the matrix-free and HBS costs for the equilibrium system $\mtx{S}$. The solve costs are expressed in terms of the actual Krylov iteration count $m$ and include unrestarted GMRES orthogonalization.}
\label{tab:S_implementation_options}
\end{table}

\begin{table}[htbp]
\centering

\begin{tabular}{@{}lrrrrrrrrr@{}}
	\toprule
	& \multicolumn{2}{c}{Memory (GB)} & \multicolumn{3}{c}{Construction (s)} & \multicolumn{2}{c}{Matvec (s)} & \multicolumn{2}{c}{Solve (s)} \\
	\cmidrule(lr){2-3}\cmidrule(lr){4-6}\cmidrule(lr){7-8}\cmidrule(lr){9-10}
	& & & & \multicolumn{2}{c}{HBS} & & & & \\
	\cmidrule(lr){5-6}
	$N$ & LU & HBS & LU & Sample & Compress & LU & HBS & LU & HBS \\
	\midrule
	$128^3$ & 14.76   & 2.38  & 12.73    & 61.48   & 30.99  & 1.33   & 0.038 & 50.50   & 1.91  \\
	$256^3$ & 314.52  & 10.32 & 366.77   & 422.92  & 61.04  & 12.16  & 0.054 & 497.79  & 4.24  \\
	$512^3$ & 6036.48 & 38.37 & 10977.61 & 3748.96 & 178.80 & 149.56 & 0.157 & 5683.28 & 15.78 \\
	\bottomrule
\end{tabular}
\caption{Matrix-free vs. HBS for 3D stencil discretizations with weak admissibility, at fixed $H=1/16$ and fixed rank $k=64$. For $N = 512^3$ the values for the LU-decomposition are extrapolated, based on the single double-wide slab statistics.}
\label{tab:lu_vs_hbs}
\end{table}

\begin{table}[htbp]\centering
\begin{tabular}{@{}lrrrrrrrrr@{}}
\toprule
& \multicolumn{2}{c}{Memory (GB)} & \multicolumn{3}{c}{Construction (s)} & \multicolumn{2}{c}{Matvec (s)} & \multicolumn{2}{c}{Solve (s)} \\
\cmidrule(lr){2-3}\cmidrule(lr){4-6}\cmidrule(lr){7-8}\cmidrule(lr){9-10}
& & & & \multicolumn{2}{c}{HBS} & & & & \\ \cmidrule(lr){5-6}
{$H$} & {LU} & {HBS} & {LU} & {Sample} & {Compress} & {LU} & {HBS} & {LU} & {HBS} \\
\midrule
{$1/8$}  & 519.96 & 4.42  & 784.87 & 414.88 & 30.66  & 14.21 & 0.023 & 261.78 & 0.87 \\
{$1/16$} & 314.52   & 10.32 & 366.77  & 422.92  & 61.04  & 12.16  & 0.054 & 497.79  & 4.24 \\
{$1/32$} & 170.32  & 22.12 & 155.59 & 505.07 & 139.21 & 9.05  & 0.111 & 635.80 & 26.84 \\
{$1/64$} & 81.39  & 45.72 & 57.81 & 520.32 & 272.16 & 8.20  & 0.239 & 1179.98 & 147.18 \\
\bottomrule
\end{tabular}
\caption{Matrix-free vs.\ HBS, for 3D stencil discretization (order $256^3$) as $H\rightarrow 0$, with fixed HBS rank $k=64$.}
\label{tab:lu_vs_hbs_Htozero}
\end{table}

The benefit of using HBS approximations for the blocks of $\mtx{S}$ is that this limits the memory footprint and enables cheap application of $\mtx{S}$. It also sets the stage for a later direct solver to be constructed (see Section~\ref{subsec:directSolve}). In contrast, the more usual (see, e.g. \cite{smithDD,toselliWidlundDD}) matrix-free approach is cheaper to set up; it requires no sampling and no compression. Note that, at a fixed number of double-wide slabs $N_{\rm ds}$, as the grid is refined, the HBS approach becomes attractive. On the other hand, if the global grid is considered fixed and the number of double-wide slabs is increased, the matrix-free approach is more interesting. However, in this case, the number of iterations required to iteratively solve the resulting equilibrium system can grow prohibitively large. For a comparison of these two regimes, see Table~\ref{tab:lu_vs_hbs} and~\ref{tab:lu_vs_hbs_Htozero}.

\subsection{Comparison to the non-overlapping approach}
\label{subsec:gmres_iters}
Section~\ref{sec:conditioning} gives an $\mathcal{O}(H^{-2})$ bound for the symmetrically scaled compatible SPD system and for the eigenvalue ratio of the raw equilibrium matrix. The supplementary material reports numerical evidence that, for representative stencil and collocation realizations, the Euclidean condition number follows a similar $H$-dependence. Motivated by these observations, the experiments below report GMRES iteration counts for the raw systems as an empirical diagnostic. We consider both SPD problems covered by the theorem and nonsymmetric or indefinite problems outside its scope. We restrict ourselves to the 3D case. All tests use a standard seven-point stencil discretization on the unit cube $\Omega=[0,1]^3$. We compare unrestarted GMRES applied to
\begin{enumerate}
    \item the original sparse system $\mtx{A}$,
    \item the nonoverlapping domain decomposition solver $\mtx{T}_H$, and
    \item the overlapping domain decomposition solver $\mtx{S}_H$.
\end{enumerate} The GMRES (relative) tolerance was set to $10^{-8}$. No $\mathcal{H}$-matrix compression was used. For the large-scale experiments, the condition number was estimated using power-iteration, rather than by explicitly computing the SVD.

Table~\ref{tab:fd_poisson_fixedH} shows the basic test for the Laplace operator $-\Delta$. At fixed $H=1/8$, mesh refinement increases the iteration counts for $\mtx{A}$ and $\mtx{T}_H$, while the $\mtx{S}_H$ iteration count remains constant. The condition number $\kappa_2(\mtx{S}_H)$ likewise remains virtually constant.
\begin{table}[H]
\centering
\begin{tabular}{c|ccc|ccc}
$N$ & $\mtx{A}$ iters & $\mtx{T}_{H}$ iters & $\mtx{S}_{H}$ iters & $\kappa_2(\mtx{A})$ & $\kappa_2(\mtx{T}_{H})$ & $\kappa_2(\mtx{S}_{H})$\\
\hline
4096 & 53 & 47 & 19 & $8.94\,10^{1}$ & $1.41\,10^{2}$ & 8.93\\
32768 & 108 & 69 & 19 & $3.84\,10^{2}$ & $2.93\,10^{2}$ & 8.89\\
262144 & 213 & 99 & 19 & $15.95\,10^{2}$ & $6.04\,10^{2}$ & 8.88\\
2097152 & 413 & 140 & 19 & $64.85\,10^{2}$ & $12.13\,10^{2}$ & 8.88
\end{tabular}
\caption{Finite-difference Poisson problem on the unit cube with eight fixed physical slabs.  The table reports unrestarted GMRES iterations for the full matrix $\mtx{A}$, the Schur-complement system $\mtx{T}_{H}$, and the scaled equilibrium system $\mtx{S}_H$.}
\label{tab:fd_poisson_fixedH}
\end{table}
In Table~\ref{tab:fd_shrinkingH} we compare the number of GMRES iterations for a decreasing slab width $H$. 
The $\mtx{S}_{H}$ iteration count grows as $H$ decreases, qualitatively consistent with the $H$-dependence of the energy and spectral bounds, while remaining smaller than the alternatives.
\begin{table}[H]
\centering
\begin{tabular}{c|ccc|ccc}
$H$ & $\mtx{A}$ iters & $\mtx{T}_{H}$ iters & $\mtx{S}_{H}$ iters & $\kappa_2(\mtx{A})$ & $\kappa_2(\mtx{T}_{H})$ & $\kappa_2(\mtx{S}_{H})$\\
\hline
1/4 & 213 & 81 & 9 & $15.96\,10^{2}$ & $6.74\,10^{2}$ & 2.45\\
1/8 & 213 & 99 & 19 & $15.96\,10^{2}$ & $7.03\,10^{2}$ & 8.88\\
1/16 & 213 & 118 & 38 & $15.96\,10^{2}$ & $8.28\,10^{2}$ & 34.82
\end{tabular}
\caption{Finite-difference Poisson problem at varying $H$. The conditioning of $\mtx{A}$ and its iteration count are included for comparison. Here $N = 64\times64\times 64$.}
\label{tab:fd_shrinkingH}
\end{table}
Table~\ref{tab:fd_other_models} provides three additional tests at fixed $H$: A variable-coefficient SPD diffusion problem, a convection-diffusion problem and a Helmholtz wavenumber sweep. 
These last two are included to examine empirically whether some of the qualitative behavior observed in the SPD tests persists in nonsymmetric and indefinite cases; they also show that, in these experiments, the growth in $\mtx{S}_{H}$-iteration count is modest as the wavenumber increases.
The Helmholtz frequency is reported in number of wavelengths $N_{\lambda}$ per unit distance, i.e. $N_{\lambda} = \kappa/(2\pi)$. Again, the underlying discretizations are $7$-point stencils in the unit cube.
\begin{table}[H]
\centering
\begin{tabular}{l|c|c|ccc}
problem & parameter & $N$ & $\mtx{A}$ iters & $\mtx{T}_H$ iters & $\mtx{S}_H$ iters\\
\hline
variable-coefficient diffusion & - & $32768$ & 136 & 79 & 19\\
variable-coefficient diffusion & - & $262144$ & 275 & 116 & 19\\
variable-coefficient diffusion & - & $2097152$ & 538 & 166 & 19\\
\hline
convection-diffusion & - & $32768$ & 122 & 70 & 20\\
convection-diffusion & - & $262144$ & 230 & 100 & 20\\
convection-diffusion & - & $2097152$ & 458 & 141 & 20\\
\hline
Helmholtz & $N_{\lambda} = 1$ & $2097152$ & $567$ & 188 & 27\\
Helmholtz & $N_{\lambda} = 2$ & $2097152$ & $915$ & 382 & 61\\
Helmholtz & $N_{\lambda} = 3$ & $2097152$ & $1000+$ & $891$ & $141$\\
\end{tabular}
\caption{Finite-difference fixed-slab diagnostics beyond constant-coefficient Poisson.  The first two blocks refine the grid at fixed physical coefficients.  The final block is a frequency sweep for Helmholtz on a $128\times 128 \times 128$ grid.}
\label{tab:fd_other_models}
\end{table}

For the variable-coefficient SPD problems, the conditioning behavior is consistent with the theory from Section~\ref{sec:conditioning}. 
The convection-diffusion and Helmholtz results provide empirical indications that similar qualitative behavior is observed in these nonsymmetric and indefinite test problems.

One problem with using a non-overlapping decomposition is that HBS compression with weak admissibility no longer suffices. In Figure~\ref{fig:gmres_comparison} we compare, for the overlapping and non-overlapping formulations $\mtx{S}_H$ and $\mtx{T}_H$, the error of the discrete Poisson BVP solution $\hat{\vct{u}}$ reached when the interface operators are compressed using weak admissibility. 

\begin{figure}
\centering
\begin{minipage}[t]{0.5\textwidth}
\centering
\vspace{0pt}
\pgfplotstableread[col sep=comma]{resultsSform_stencil.csv}\datatableS
\pgfplotstableread[col sep=comma]{resultsTform_stencil.csv}\datatableT
\begin{tikzpicture}
\begin{axis}[
    width=\linewidth,
    height=0.5\linewidth,
    xlabel={$H$},
    xlabel shift=-7pt,
    ylabel={$\|\vct{u}-\widehat{\vct{u}}\|_2/\|\vct{u}\|_2$},
    ymode=log,
    xmode=log,
    log basis x={2},
    grid=both,
]
\addplot[
    black,
    thick,
    mark=o,
]
table[
    x=H,
    y=gmres_err
]{\datatableS};
\addplot[
    black,
    thick,
    mark=square*,
]
table[
    x=H,
    y=gmres_err
]{\datatableT};
\end{axis}
\end{tikzpicture}

\pgfplotstableread[col sep=comma]{rkErrS.csv}\datatableSrk
\pgfplotstableread[col sep=comma]{rkErrT.csv}\datatableTrk
\begin{tikzpicture}
\begin{axis}[
    width=\linewidth,
    height=0.5\linewidth,
    xlabel={rank},
    ylabel={$\|\vct{u}-\widehat{\vct{u}}\|_2/\|\vct{u}\|_2$},
    ymode=log,
    grid=both,
    legend pos = north east
]
\addplot[
    black,
    thick,
    mark=o,
]
table[
    x=rk,
    y=err
]{\datatableSrk};
\addlegendentry{$S$}
\addplot[
    black,
    thick,
    mark=square*,
]
table[
    x=rk,
    y=err
]{\datatableTrk};
\addlegendentry{$T$}
\end{axis}
\end{tikzpicture}
\end{minipage}
\hfill
\begin{minipage}[t]{0.45\textwidth}
\small
\vspace{0pt}
\textit{Accuracy of the GMRES solution obtained by $\mtx{T}_{H}$ and $\mtx{S}_H$ for a $128\times128\times128$ stencil discretization of a Poisson problem in $\Omega = [0,1]^3$. Top: fixed rank set to $128$, slab spacing $H$ varying. Bottom: fixed slab width $H=1/8$, HBS rank varying. The operator $\mtx{T}_H$, as it contains pseudo-differential self-interactions of the interfaces, requires prohibitively high ranks to converge. In this experiment, the overlapping formulation exhibits substantially milder conditioning and iteration behavior than the non-overlapping formulation. The $\mtx{S}_{H}$-error also increases as $H\rightarrow 0$, but remains much more benign.}
\end{minipage}
\caption{Comparison of solution error for formulations $\mtx{S}$ and $\mtx{T}$
as a function of slab width $H$ and the HBS compression rank.}
\label{fig:gmres_comparison}
\end{figure}
In summary, the reported experiments show that the non-overlapping approach not only has larger iteration count, but also shows much poorer HBS-compressibility (with weak admissibility) than its overlapping counterpart.
\begin{remark}
The fact that the overlapping formulation can ``get away'' with using weak admissibility in 3D (as shown in Figure~\ref{fig:gmres_comparison}) is a particular feature of our method, and results from the fact that we deliberately sought a formulation that involves integral operators with \textit{smooth} kernels, cf.~(\ref{eq:debbie1}) and (\ref{eq:debbie2}). However, other rank structured formats can easily be used; either simpler single level structures \cite{2019_amestoy_BLR_multifrontal}, or more complex ones such as $\mathcal{H}^2$-matrices with strong admissibility \cite{2002_hackbusch_H2,2010_borm_book}.
\end{remark}
\subsection{Direct solvers}\label{subsec:directSolve}
As the above demonstrates, the use of iterative solvers can be disadvantageous whenever $H\rightarrow 0$ or the problem becomes highly oscillatory. The benefit of the reduced $\mtx{S}$-system then is that, in the absence of (numerical) interior resonances, it is a well-posed block tridiagonal system for which various solution strategies exist. These are based on combining block tridiagonal solvers (e.g., the Thomas algorithm or cyclic reduction) with $\mathcal{H}$-matrix arithmetic. This will be the subject of future work.
\section{Numerical examples}\label{sec:numerical}
In this section we demonstrate the effectiveness of the proposed method on some challenging 2D and 3D examples. The experiments were performed on a 24 core Intel Xeon Gold 6248R 3GHz CPU machine, using our Python software package SslabLU (available at~\cite{SslabLU}).
The HBS construction timings we will report in Example~\ref{sub:example_2_photonic_crystal_waveguide} and~\ref{subsec:twisted} are broken down into:
\begin{enumerate}
    \item $t_{\mtx{A}}$, the total time to construct and factorize all local stiffness matrices
    \item $t_{\rm sample}$, the total time over the double-wide slabs to probe $k$ random samples needed in the construction of the HBS approximants to $\mtx{S}_{j,j\pm1}$
    \item $t_{\text{HBS}}$, the total time to compress the blocks $\mtx{S}_{j,j'}$ over the double-wide slabs after sampling has been performed.
\end{enumerate}
We note that the number of double slabs $N_{\rm ds}$ corresponds to the number of interfaces. It is equal to the number of single slabs in the periodic case, and to one less than the number of single slabs in the non-periodic case. We also stress that the timings are for an unparallelized implementation. As outlined in Remark~\ref{rem:parallel}, each of these three steps can be parallelized.

\subsection{Example 1: Variable-coefficient screened diffusion equations}\label{sec:convdiff}
Screened diffusion equations of the form
\begin{align}\label{eq:screened_diff}
-\nabla\cdot(D(x,y)\nabla u(x,y)) + \alpha^2u(x,y) &= f(x,y) &\text{ in } \Omega\\
\frac{\partial u(x,y)}{\partial n} &= 0 &\text{ on }\partial\Omega
\end{align}
are a powerful tool in the numerical solution of convection-diffusion equations in incompressible flow:
\begin{align*}
    \frac{\partial u(t;x,y)}{\partial t} + \vct{v}(t;x,y)\cdot\nabla u(t;x,y) &= \nabla\cdot(D(x,y)\nabla u(t;x,y)) & \text{ in }[0,t_{\rm fin}]\times\Omega\\
    \frac{\partial u(t;x,y)}{\partial n} &= 0 & \text{ on }\partial\Omega\\
    u(0;x,y)&=U_0(x,y) & \text{ in } \quad \Omega.
\end{align*}
By introducing an implicit-explicit (IMEX) time-stepping discretization, as in~\cite{FRYKLUND2023111856}, \eqref{eq:convdiff} can be reduced to a sequence of problems of the form \eqref{eq:screened_diff}. To be explicit:
\begin{align}\label{eq:convdiff}
-\nabla\cdot(D(x,y)\nabla U_{n+1}(x,y)) + \frac{1}{\Delta t}U_{n+1}(x,y) &= \frac{1}{\Delta t}U_n(x,y)-\vct{v}\cdot\nabla U_n(x,y) \text{ in } \Omega
\end{align}
with zero Neumann boundary conditions and $U_{n}(x,y) = u(n\Delta t;x,y)$. This is known as \textit{elliptic marching}, \textit{Rothe's method} or the \textit{method of lines}. For simplicity, we restrict our attention here to a first order IMEX scheme (backward Euler), but more advanced time-stepping schemes and IMEX discretizations can be used. See~\cite{FRYKLUND2023111856} and ~\cite{ARWpaper} for details on spectral correction time-stepping and IMEX methods.

We take $\Omega = [0,1]^2$ with periodic boundary conditions in the horizontal direction, representing an infinitely long horizontal channel. We set $D(x,y) = D_0(1 + \varepsilon(s_y(2y-1) + s_x \cos(2\pi x)))$ with $D_0=0.1,\varepsilon=0.3,s_y=1 \text{ and } s_x=0.5$, and $\vct{v}(x,y)=(-2y(1-y),0)$. The initial concentration $U_0$ is set to be a smooth bump function. 

The system in \eqref{eq:convdiff} was discretized using an HPS discretization consisting of $16 \times 16$ square subdomains, carrying a $p\times p$ Chebyshev collocation grid, for $p\in\{6,\ldots,16\}$. We study the self-convergence of $u(1;x,y)$ at $t_{\rm fin}=1$, with respect to a high order reference solution at $p_{\rm ref}=20$. We take $\Delta t = 6.94\times 10^{-3}$.\footnote{This is the timestep needed to satisfy the CFL condition at the reference polynomial order $p_{\rm ref}=20$.} Since $\Delta t$ is positive and fixed, the continuum operator on the left-hand side of \eqref{eq:convdiff} is always symmetric positive definite. The HPS collocation matrix is in general not SPD in the Euclidean inner product, but the same $\mtx{S}$-matrix reduction can be reused for all time-steps. Slab separation $H = 1/16$ was used for each order. The reference solution was obtained using the matrix-free approach, i.e., without any compression.

\begin{figure}[H]
    \centering
    \includegraphics[width=0.32\linewidth]{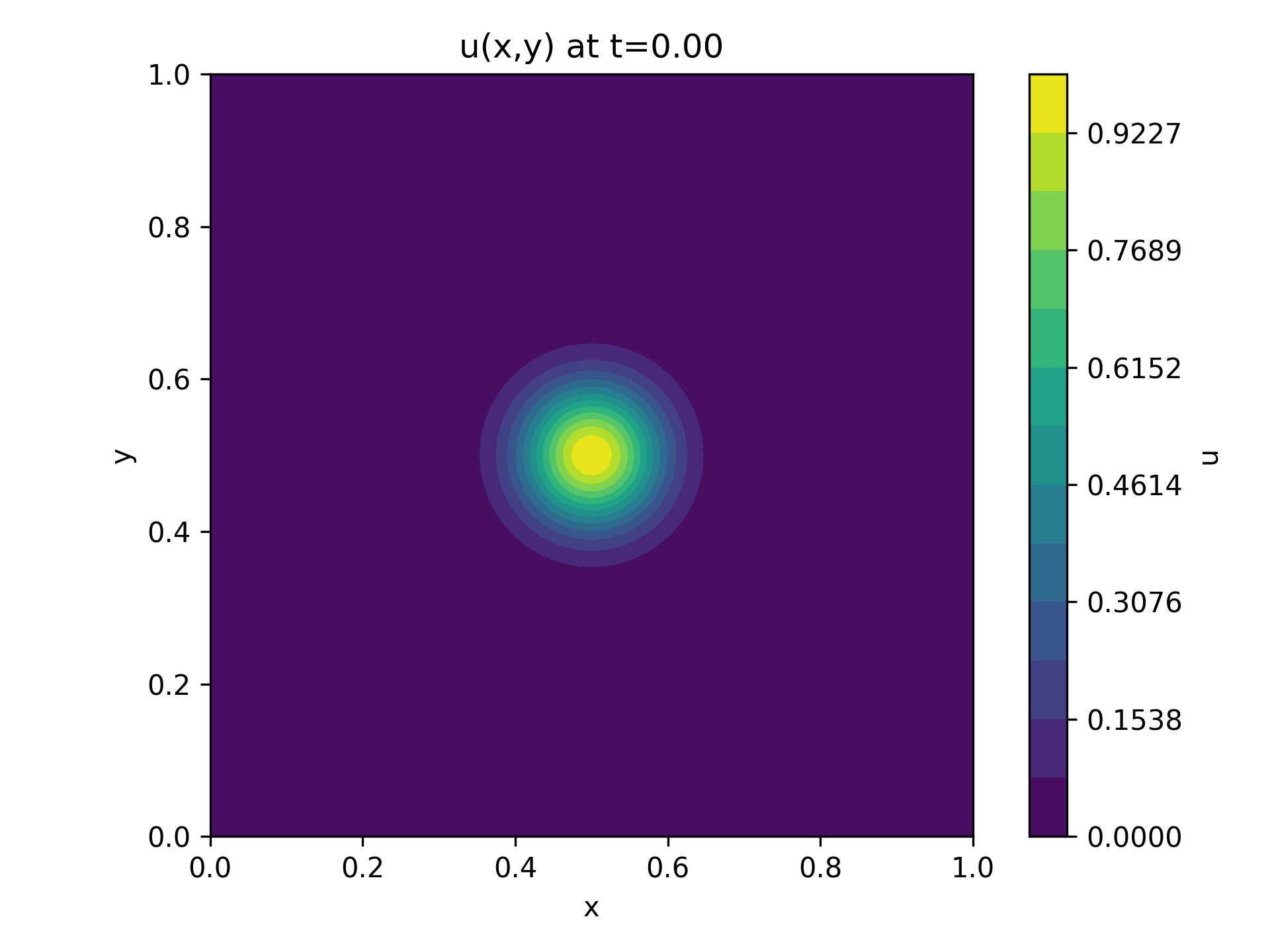}
    \hfill
    \includegraphics[width=0.32\linewidth]{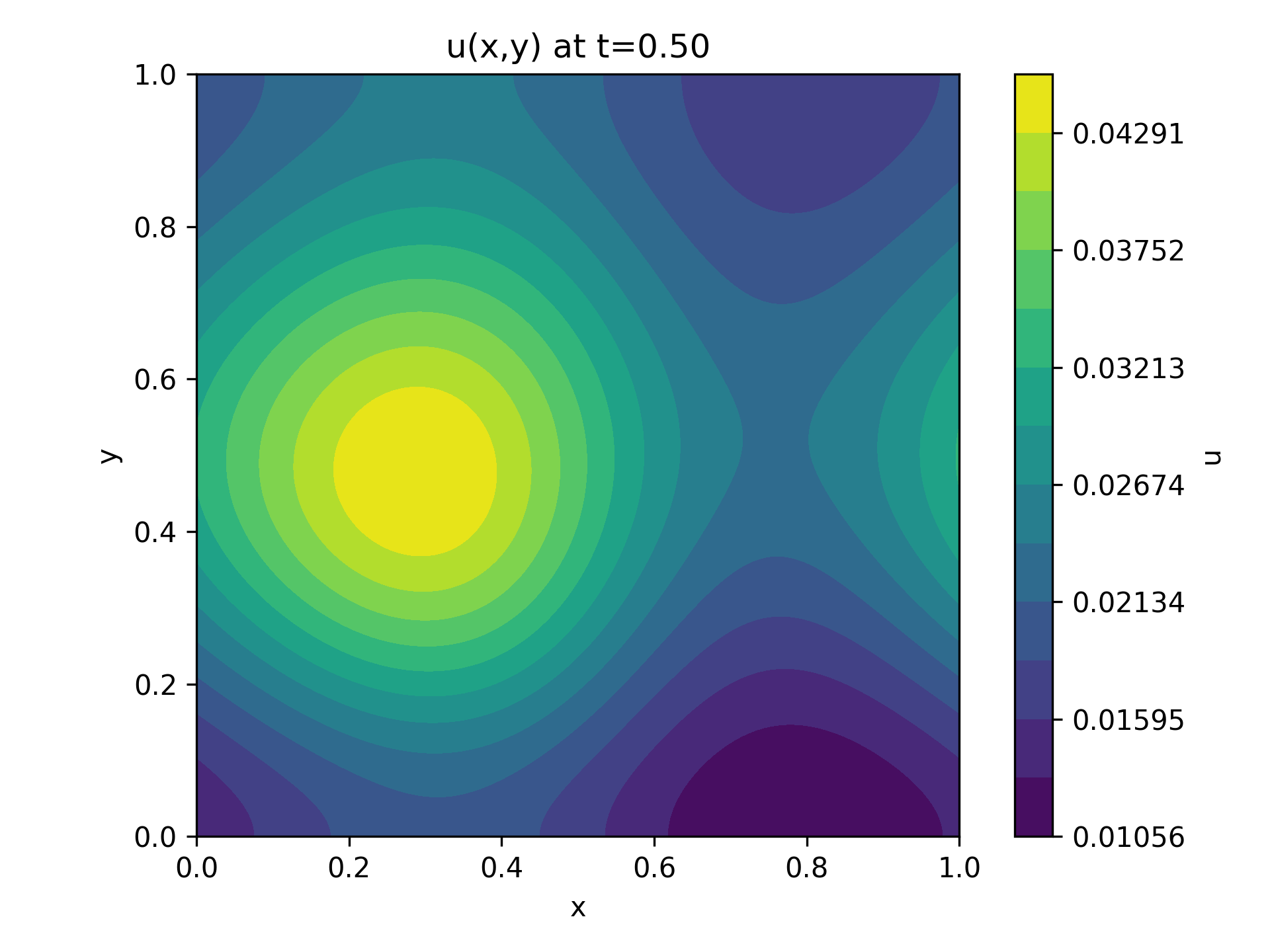}
    \hfill
    \includegraphics[width=0.32\linewidth]{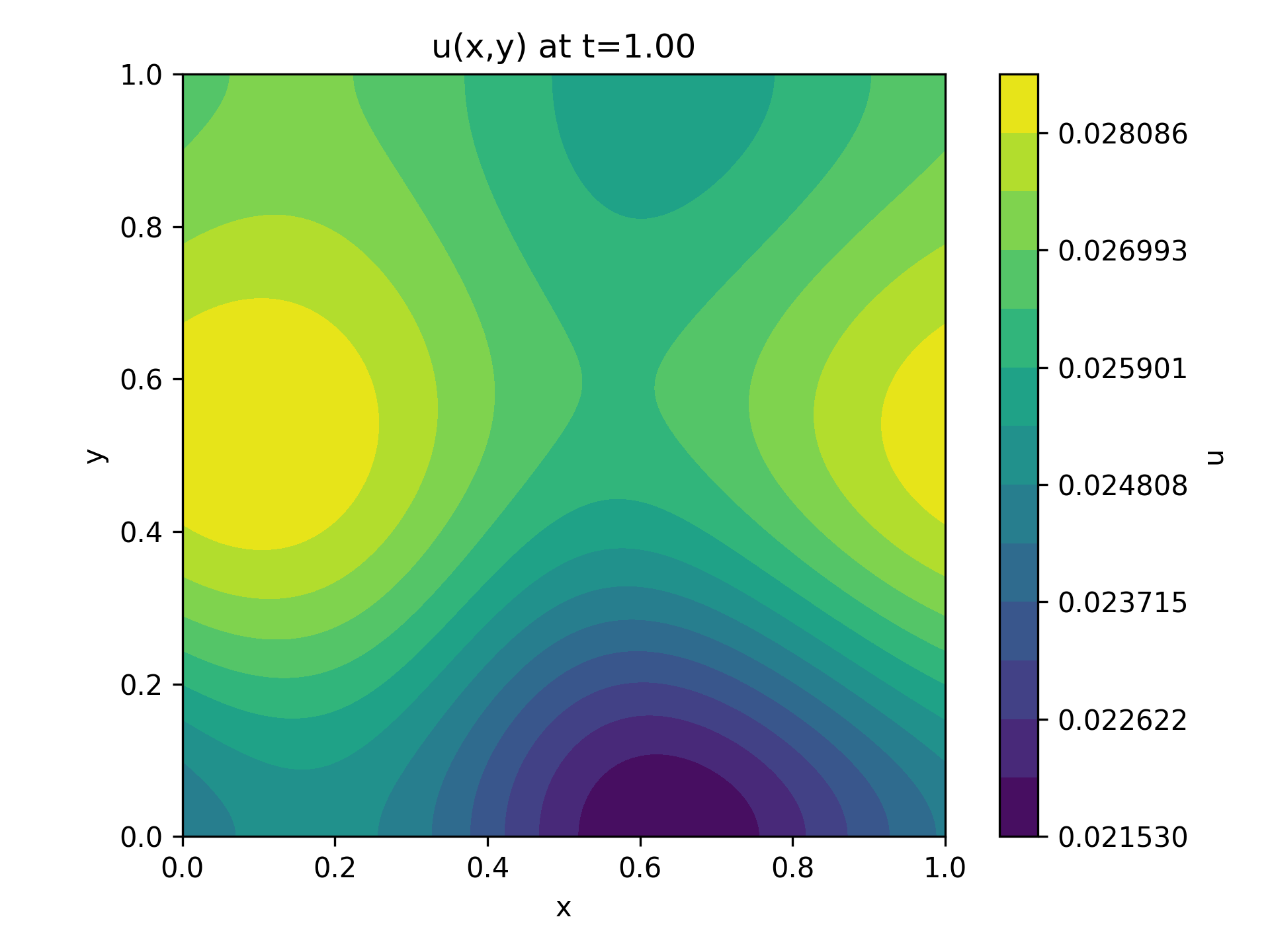}
    \caption{Contour plots of the solution $U_n$ to \eqref{eq:convdiff} at $t=0$, $t=.5$ and $t=1.$}
    \label{fig:convdiffPlots}
\end{figure}

Figure~\ref{fig:convdiff_comparison} reports, at various $N=16^2p^2$, the error of the discrete solution interpolated on a fine equispaced grid with respect to a reference solution interpolated on that same grid. It also shows the GMRES solve time for the equilibrium system $\mtx{S}$ averaged over the IMEX time-steps. In Figure~\ref{fig:convdiffPlots} the solution $U_{n}(x,y)$ (approximated by the reference solution) is plotted for $n\Delta t\in\{0\,,\,1/2\,,\,1\}$. The timing results for the construction of the HBS approximation are completely similar to those in Section~\ref{sub:example_2_photonic_crystal_waveguide}, hence they are not included here.

\begin{figure}
\centering
\begin{minipage}[t]{0.55\textwidth}
\centering
\vspace{0pt}
\pgfplotstableread[col sep=comma]{dataOMS/convdiffstats.csv}\datatableConvDiff
\begin{tikzpicture}
\begin{groupplot}[
    group style={
        group name=my plots,
        group size=1 by 2,
        xlabels at=edge bottom,
        xticklabels at=edge bottom,
    },
    footnotesize,
    width=\linewidth,
    height=4cm,
    xlabel=$N$,
    ytick align=outside,
    grid=both,
    xmode=log,
    ymode=log,
    legend pos = north west
]
\nextgroupplot[title={Error at $t=1$ w.r.t. reference solution},ylabel=$\|\vct{u}_{\rm ref}-\vct{u}_p\|/\|\vct{u}_{\rm ref}\|$]
\addplot[
    black,
    thick,
    mark=o,
]
table[
    x=N,
    y=err_HBS
]{\datatableConvDiff};
\addplot[
    black,
    thick,
    mark=square,
]
table[
    x=N,
    y=err_lu
]{\datatableConvDiff};
\nextgroupplot[title={Avg. solve time per IMEX iteration},ylabel=time (s)]
\addplot[
    black,
    thick,
    mark=o,
]
table[
    x=N,
    y=t_solve_HBS
]{\datatableConvDiff};
\addlegendentry{HBS}
\addplot[
    black,
    thick,
    mark=square,
]
table[
    x=N,
    y=t_solve_LU
]{\datatableConvDiff};
\addlegendentry{Matrix-free}
\end{groupplot}
\end{tikzpicture}
\end{minipage}
\hfill
\begin{minipage}[t]{0.38\textwidth}
\small
\vspace{10pt}
\textit{The convection-diffusion equation \eqref{eq:convdiff} is solved for $n=1,\ldots,n_{\rm fin}$ such that $n_{\rm fin}\Delta t = 1$. The vector $\vct{u}_p$ denotes the approximation of $U_{n_{\rm fin}}(x,y)$ using a $16\times 16$ HPS grid of order $p$ at every IMEX time-step. We track the average GMRES solution time per IMEX time-step, as well as the final error at $n_{\rm fin}$ with respect to a reference solution $\vct{u}_{\rm ref}=\vct{u}_{20}$. For both the matrix-free formulation and the HBS compressed formulation, the number of GMRES iterations at each timestep was (at most) $9$. The HBS rank was set to $k=16$ for each $N=16^2p^2$. At this rank, the HBS compression is essentially exact. This is reflected in the top figure; the two error graphs cannot be distinguished.}
\end{minipage}
\caption{Error at $t=1$ and timing of the matrix-free and HBS-compressed $\mtx{S}$-formulation for $u(1,x,y)$.}
\label{fig:convdiff_comparison}
\end{figure}
\subsection{Example 2: Photonic crystal waveguide}\label{sub:example_2_photonic_crystal_waveguide}

\begin{figure}
    \centering
    \begin{subfigure}{.32\linewidth}
        \centering
        \includegraphics[width=\linewidth]{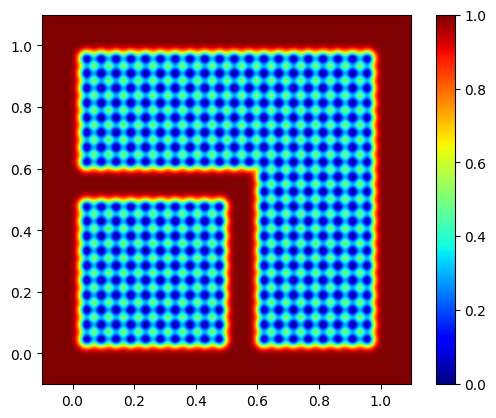}
        \caption{$1-b(x,y)$}
        \label{subfig:bfield}
    \end{subfigure}
    \hfill
    \begin{subfigure}{.32\linewidth}
    \centering
        \includegraphics[width=\linewidth]{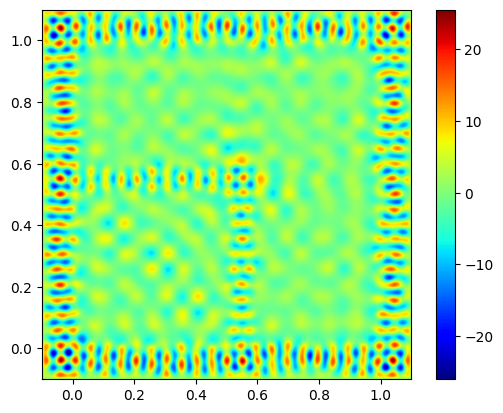}
        \caption{Solution at $\kappa = 156.703$}
    \end{subfigure}
    \begin{subfigure}{.32\linewidth}
    \centering
        \includegraphics[width=\linewidth]{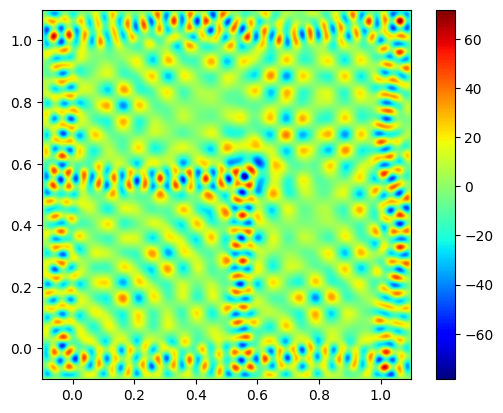}
        \caption{Solution at $\kappa = 157.017$}
        \label{subfig:waveguide-kappa-157}
    \end{subfigure}
    \caption{Numerical solution, at $\kappa = 156.703$ and $\kappa = 157.017$, to the BVP in equation~(\ref{eq:waveGuide}).}
    \label{fig:enter-label}
\end{figure}

We solve the Helmholtz equation with variable wavenumber:
\begin{equation}\label{eq:waveGuide}
\begin{cases}
    -\Delta u(x,y)-\kappa^2(1-b(x,y))u(x,y) = 0, & (x,y)\in\Omega,\\
    u(x,y)=1, & (x,y)\in\partial\Omega,
\end{cases}
\end{equation}
where $\Omega=[0,1]^2$ and $1-b$ models the relative speed of light imposed by a crystal waveguide, represented as a collection of Gaussian bumps (see Figure~\ref{subfig:bfield}).

We set $H=1/8$ and the HPS tiling per double-wide slab to be an $8$-by-$32$ square lattice and study \textit{self-convergence} of the solution: for $p\in\{8,10,\ldots,18,20\}$ we compare, on a fine uniform grid, the interpolated solution $\vct{u}_{p}$ to an interpolated reference solution $\vct{u}_{30}$ at $p=30$ (computed without HBS compression and using a direct solver). We use HBS compression with rank $k=25$. We use GMRES as an iterative solver, with its tolerance set to $10^{-10}$. For reference, we also plot the self-convergence \textit{without} HBS compression and using a direct solver. In a sense this provides the best convergence we could hope for.
\begin{figure}[H]
    \centering
    \begin{subfigure}[t]{.46\linewidth}
    \begin{tikzpicture}
    \begin{axis}[
    width=\linewidth,
    height=6cm,
    xmode = log,
    ymode = log,
    grid style=dashed,
    grid =both,
    ylabel={$\|\vct{u}_{p}-\vct{u}_{30}\|/\|\vct{u}_{30}\|$},
    xlabel={$N$},
    ytick={1e-10,1e-8,1e-6,1e-4,1e-2,1e-0},
    xtick={1024*(8^2),1024*(12^2),1024*(16^2),1024*(20^2)},
    legend pos = north east
]
\addplot[black, mark=*
    ] table [x expr = 1024*\thisrow{p}^2,y=err,col sep=comma]{dataOMS/crystal_waveguide_hmat.csv};
\addplot[black, mark=*,dashed
    ] table [x expr = 1024*\thisrow{p}^2,y=err,col sep=comma]{dataOMS/crystal_waveguide_dense.csv};
    \legend{HBS,dense}
    \end{axis}
    \node[above] at (.45,4.5){\footnotesize$p=8$};
    \node[above] at (2.3,4.5){\footnotesize$p=12$};
    \node[above] at (3.7,4.5){\footnotesize$p=16$};
    \node[above] at (4.8,4.5){\footnotesize$p=20$};
    \end{tikzpicture}    
    \caption{Self-convergence error $\|\vct{u}_{p}-\vct{u}_{30}\|/\|\vct{u}_{30}\|$ as a function of $p$.}
    \end{subfigure}
\hfill
\begin{subfigure}[t]{.46\linewidth}
\pgfplotsset{select coords between index/.style 2 args={
    x filter/.code={
        \ifnum\coordindex<#1\def\pgfmathresult{}\fi
        \ifnum\coordindex>#2\def\pgfmathresult{}\fi
    }
}}
    \begin{tikzpicture}
    \begin{axis}[
    width=\linewidth,
    height=6cm,
    ymin=.3,
    ylabel={time $[s]$},
    xlabel={$N$},
    ymode=log,
    xmode=log,
    log basis x={10},
    grid = both,
    grid style=dashed,
    legend pos = north west,
    xtick={1024*(8^2),1024*(12^2),1024*(16^2),1024*(20^2)},
]
\addplot[black, mark=square*
    ] table [x expr =\thisrow{p}*\thisrow{p}*1024,y expr=7*\thisrow{discr},col sep=comma]{dataOMS/crystal_waveguide_hmat.csv};
\addplot[black, mark=diamond*
    ] table [x expr=\thisrow{p}*\thisrow{p}*1024,y expr = 7*\thisrow{sample},col sep=comma]{dataOMS/crystal_waveguide_hmat.csv};
\addplot[black, mark=*
    ] table [x expr=\thisrow{p}*\thisrow{p}*1024,y expr = 7*\thisrow{compr},col sep=comma]{dataOMS/crystal_waveguide_hmat.csv};
\addplot[black,dashed,select coords between index={4}{6}
    ] table [x expr=\thisrow{p}*\thisrow{p}*1024,y expr=7*\thisrow{p}/180,col sep=comma]{dataOMS/crystal_waveguide_hmat.csv}node[below,pos=.5,rotate=11] {$\mathcal{O}(N^{1/2})$};
    \addplot[black,dashed,select coords between index={4}{6}
    ] table [x expr=\thisrow{p}*\thisrow{p}*1024,y expr=7*\thisrow{p}^2*log2(\thisrow{p}^2)/10000,col sep=comma]{dataOMS/crystal_waveguide_hmat.csv}node[above,pos=.5,rotate=30] {$\mathcal{O}(N\log N)$};
    \addplot[black,dashed,select coords between index={2}{6}
    ] table [x expr=\thisrow{p}*\thisrow{p}*1024,y expr=7*\thisrow{p}^3/3000,col sep=comma]{dataOMS/crystal_waveguide_hmat.csv}node[above,pos=.25,rotate=30] {$\mathcal{O}(N^{3/2})$};
    \legend{$t_{\mtx{A}}$,$t_{\rm sample}$,$t_{\text{HBS}}$}
    \end{axis}
    
    \node[above] at (.45,4.5){\footnotesize$p=8$};
    \node[above] at (2.3,4.5){\footnotesize$p=12$};
    \node[above] at (3.7,4.5){\footnotesize$p=16$};
    \node[above] at (4.8,4.5){\footnotesize$p=20$};
    \end{tikzpicture}    
    \caption{Timings $t_{\mtx{A}}$, $t_{\rm sample}$ and $t_{\text{HBS}}$ as a function of the total DOFs $N$. Expected asymptotic costs are shown with dashed lines.}
    \end{subfigure}    
    \caption{Convergence and computational cost for the solution of equation~(\ref{eq:waveGuide}), using our proposed solver. Here $N_{\rm ds} = 7$ and each double-wide slab uses an $8$-by-$32$ HPS tiling. The polynomial order $p$ is indicated on top. The accuracy stalls around $6-7$ digits because the off-diagonal block compression is accurate up to 7 digits. Due to ill-conditioning of the internal solver, convergence for the dense solver slows down at high $p$.}
    \label{fig:waveGuideStats}
\end{figure}
We observe the following:
\begin{enumerate}
    \item For all $p$ the number of GMRES iterations required for both wave numbers was $198 \pm 1$. Thus, for this fixed physical problem and slab geometry, the GMRES iteration count is essentially unchanged as the HPS order $p$ varies.
    \item Even with HBS compression at rank `only' $k=25$, we see that the convergence is essentially optimal, stalling around $10^{-5}$, the accuracy threshold at this rank. For this 2D problem then, despite its oscillatory nature, HBS compression does not influence convergence, only the maximal accuracy reached.
\end{enumerate}
HBS compression, including sampling, is negligible compared to building and factorizing the local stiffness matrices; both factorization and sampling costs match the predicted asymptotics.
\begin{remark}
    A common approach to solving the Helmholtz equation at high frequency is to construct Impedance-to-impedance maps instead of Dirichlet-to-Neumann or our Dirichlet-to-Dirichlet maps. In Appendix~\ref{app:example2_conditioning} we investigate the GMRES iteration counts for these three approaches.
\end{remark}

\subsection{Example 3: Twisted square torus}\label{subsec:twisted}
Another source of variable coefficient PDEs is transformed constant coefficient PDEs. Say we wish to solve the 3D Helmholtz equation
\begin{align}\label{eq:HelmholtzTorus}
    -\Delta u(x,y,z) - \kappa^2 u(x,y,z) &= 0\text{ in }\Omega\\
    u(x,y,z) &= 1\text{ on }\partial\Omega
\end{align}
where $\Omega$ admits transformations
\begin{align*}
    f&:[0,1]^3\rightarrow\Omega\\
    g&:\Omega\rightarrow[0,1]^3
\end{align*}
such that $g \circ f=\mathit{Id}$ and $f$ and $g$ are twice continuously differentiable. We consider here $\Omega$ the \textit{twisted torus}, shown in Figure~\ref{fig:twistedTorus}. Because the twisted torus is periodic, we have that $g$ and $f$ are periodic. The Helmholtz equation on $\Omega$ is transformed into a (periodic) variable coefficient PDE on $[0,1]^3$ using $g$ and $f$. We solve equation~(\ref{eq:HelmholtzTorus}) at $\kappa = 17.66$. For the untransformed twisted square torus this makes the original \textit{dimensionless wavenumber} $\kappa\cdot\text{diam}(\Omega)\approx135.46$. This means $\Omega$ is $21.56$ wavelengths across.
We again study self-convergence for $p\in\{4,6,8,10\}$, comparing to a reference solution $\vct{u}_{12}$. In all cases we use an HPS discretization in the double slabs $\Psi_j$ with an $8\times16\times16$ HPS tiling. Observe that, for $p=12$ each slab has a total number of $3,538,944$ DOFs before reduction. We set $H=1/16$ meaning $\Omega$ is discretized using $N=28,311,552$ total DOFs. For the approximate solutions $\vct{u}_{p}$ the HBS rank was set to $k=150$, and the GMRES convergence tolerance was set to $10^{-6}$. 

We observe the following:
\begin{enumerate}
    \item For all $p$ the number of GMRES iterations required was $1476$. This again demonstrates that the number of iterations is independent of the order $p$ but it is still quite high; it motivates the development of an efficient direct solver or of a cheap preconditioner.
    \item Again, we see fast convergence of the solution up to the GMRES tolerance and the accuracy of the HBS block compression. 
\end{enumerate}
We also see that the predicted asymptotic costs are slight overestimates of the observed computational costs.
\begin{figure}
    \centering
    \begin{subfigure}[t]{.42\linewidth}
    \centering
    \raisebox{.25\height}{
    \includegraphics[width=.9\linewidth]{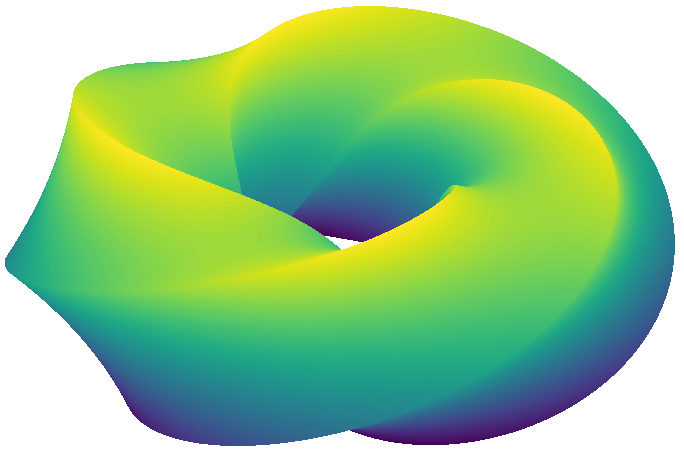}
    }\caption{Boundary of the twisted square torus domain in $\mathbb{R}^3$.}
    \label{fig:twistedTorus}
    \end{subfigure}
    \hfill
    \begin{subfigure}[t]{.42\linewidth}
    \centering
    \includegraphics[width=\linewidth]{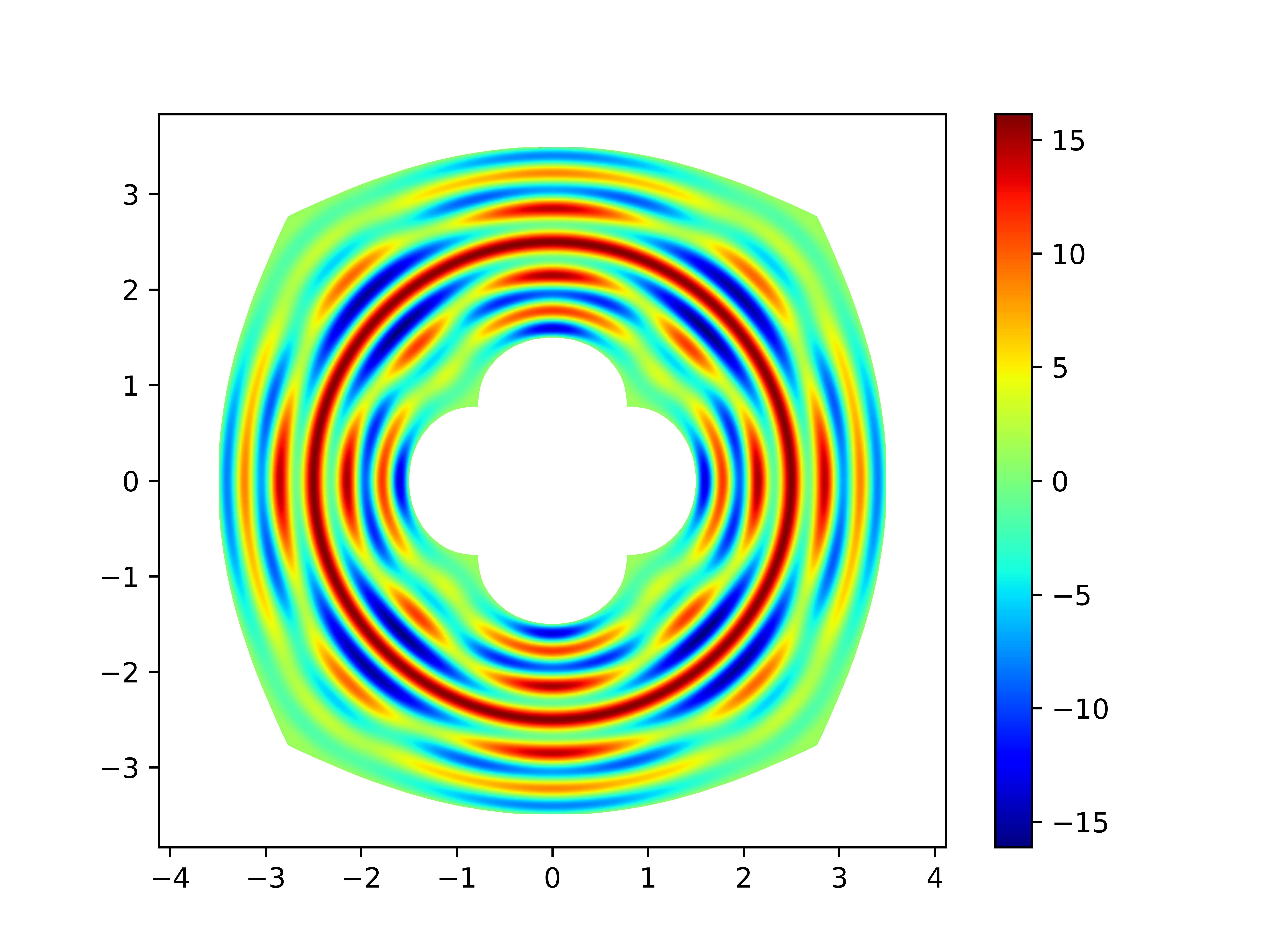}
    \caption{Cross section at $z=0$ of the solution at $\kappa=17.66$ for the Helmholtz equation on the twisted torus.}
    \end{subfigure}
    \caption{Boundary of the twisted square torus (left) and cross section of the solution at $\kappa=17.66$ of the Helmholtz equation~(\ref{eq:HelmholtzTorus}), obtained using an HPS discretization of order $p=12$.}
    \label{fig:twistedTorus}
\end{figure}
\begin{figure}
    \centering
    \begin{subfigure}[t]{.46\linewidth}
    \begin{tikzpicture}
    \begin{axis}[
    width=\linewidth,
    xmode = log,
    ymode = log,
    minor grid style={draw=none},
    xtick = {16384*(4^3),16384*(6^3),16384*(8^3),16384*(10^3)},
    grid style=dashed,
    grid = both,
    ylabel={$\|\vct{u}_{p}-\vct{u}_{12}\|/\|\vct{u}_{12}\|$},
    xlabel={$N$},
    legend pos = south east
]
\addplot[black, mark=*
    ] table [x expr = 16384*\thisrow{p}^3,y=err,col sep=comma]{dataOMS/twistedTorus.csv};
    \end{axis}
    \node[above] at (.45,4.5){\footnotesize$p=4$};
    \node[above] at (2.3,4.5){\footnotesize$p=6$};
    \node[above] at (3.7,4.5){\footnotesize$p=8$};
    \node[above] at (4.8,4.5){\footnotesize$p=10$};
    \end{tikzpicture}    
    \caption{Self-convergence error $\|\vct{u}_{p}-\vct{u}_{12}\|/\|\vct{u}_{12}\|$ as a function of $p$.}
    \end{subfigure}
\hfill
\begin{subfigure}[t]{.46\linewidth}
    \begin{tikzpicture}
    \pgfplotsset{select coords between index/.style 2 args={
    x filter/.code={
        \ifnum\coordindex<#1\def\pgfmathresult{}\fi
        \ifnum\coordindex>#2\def\pgfmathresult{}\fi
    }
}}
    \begin{axis}[
    width=\linewidth,
    ylabel={time $[s]$},
    xlabel={$N$},
    xtick = {16384*(4^3),16384*(6^3),16384*(8^3),16384*(10^3)},
    ymode=log,
    xmode = log,
    grid style=dashed,
    grid = both,
    legend pos = north west
]
\addplot[black, mark=*
    ] table [x expr = (\thisrow{p}^3)*16384,y expr = 8*\thisrow{discr},col sep=comma]{dataOMS/twistedTorus.csv};
    \addplot[black, mark=square*
    ] table [x expr = (\thisrow{p}^3)*16384,y expr = 8*\thisrow{sample},col sep=comma]{dataOMS/twistedTorus.csv};
    \addplot[black, mark=square*
    ] table [x expr = (\thisrow{p}^3)*16384,y expr = 8*\thisrow{compr},col sep=comma]{dataOMS/twistedTorus.csv};
    \addplot[black,dashed
    ] table [select coords between index={1}{5},x expr = (\thisrow{p}^3)*16384,y expr = 8*(\thisrow{p}^4)/75,col sep=comma]{dataOMS/twistedTorus.csv}node[pos=.5] (OH) {};
    \node [below,rotate=25] at (OH) {$\mathcal{O}(N^{4/3})$};
    \addplot[black,dashed
    ] table [select coords between index={1}{5},x expr = (\thisrow{p}^3)*16384,y expr = (\thisrow{p}^6)/75,col sep=comma]{dataOMS/twistedTorus.csv}node[pos=.5] (OH2) {};
    \node [above,rotate=30] at (OH2) {$\mathcal{O}(N^{2})$};
    \addplot[black,dashed
    ] table [select coords between index={1}{5},x expr = (\thisrow{p}^3)*16384,y expr = 8*(\thisrow{p}^2)/15,col sep=comma]{dataOMS/twistedTorus.csv}node[pos=.7] (OH2) {};
    \node [below,rotate=10] at (OH2) {$\mathcal{O}(N^{2/3})$};
    \legend{$t_{\mtx{A}}$,$t_{\rm sample}$,$t_{\text{HBS}}$}
    \end{axis}
    \node[above] at (.45,4.5){\footnotesize$p=4$};
    \node[above] at (2.3,4.5){\footnotesize$p=6$};
    \node[above] at (3.7,4.5){\footnotesize$p=8$};
    \node[above] at (4.8,4.5){\footnotesize$p=10$};
    \end{tikzpicture}    
    \caption{Timings $t_{\mtx{A}}$, $t_{\rm sample}$ and $t_{\text{HBS}}$ as a function of the total DOFs $N$. Expected asymptotic costs are shown with dashed lines.}
    \end{subfigure}    
    \caption{Convergence and computational cost for the solution of equation~(\ref{eq:HelmholtzTorus}), using our proposed solver, for $\kappa = 17.66$. Here $N_{\rm ds} = 16$ and each double-wide slab uses an $8$-by-$16$-by-$16$ HPS tiling. The polynomial order $p$ is indicated on top. We see that convergence is slowed down as the error reaches the desired GMRES tolerance and the accuracy of the block compression.}
    \label{fig:waveGuideStats}
\end{figure}

\section{Conclusions and future work}
This manuscript describes a technique for solving scalar elliptic PDEs, based on an overlapping domain decomposition method built from thin slabs.
The method is based on an interface equilibrium equation whose off-diagonal blocks discretize interface solution operators between separated interfaces.  This structure explains the strong compressibility observed in the experiments and the favorable behavior in the numerical tests. The non-zero blocks in the equilibrium system approximate integral operators with \textit{smooth} kernels, making them highly amenable to compression using hierarchical matrix techniques.
In our method, these blocks are formed using a randomized compression technique coupled with a local direct solver that exploits that each subdomain is thin.
The numerical experiments demonstrate that the formulation can be applied accurately at large scale, including to oscillatory test problems, while also indicating the regimes where additional preconditioning or direct solvers are needed.

For coercive problems, multigrid methods remain the natural benchmark for one-off low-order full-volume solves.  The present formulation has a different cost profile: it is high-order, reusable across right-hand sides, naturally parallel across slabs, and can compute interface traces before any volume reconstruction. These features make time-to-accuracy and multiple-right-hand-side comparisons an important topic for future work.

In this work, the reduced global system is solved using an iterative method that typically converges rapidly.
However, it is also possible to construct a linear complexity fully direct solver by exploiting the rank structure in the system to compute an LU factorization in ``data sparse'' form; work in this direction is in progress.
Additionally, the technique is being extended to different boundary conditions (e.g., mixed, impedance,\ldots), and these options will be added to our software as they become available.
We are also working on developing an HPC distributed memory implementation of the solver, leveraging that the method we present is highly parallelizeable.

\section{Acknowledgements}
The work reported was supported by the Office of Naval Research (N00014-18-1-2354), by the National Science Foundation (DMS-2313434), by the Department of Energy ASCR (DE-SC0025312), and by the Sloan Foundation.
The authors would also like to express their gratitude to Joseph Kump for his contributions to the software package \texttt{SslabLU}.

\clearpage
\bibliographystyle{abbrv}
\bibliography{main_bib}
\clearpage
\appendix
\section{Trace-system iteration comparison for Example 2}\label{app:example2_conditioning}
We compare three reduced trace systems for the example in \ref{sub:example_2_photonic_crystal_waveguide}, with $\kappa=157.017$ and local HPS order $p=8$.  All three systems use the same physical slab interfaces for domain decomposition with $N_{\rm ds} = 7$.  The physical solution for this case is shown in Figure~\ref{subfig:waveguide-kappa-157}.  The first system is the standard non-overlapping DtN Schur complement $\mtx{T}$ from (\ref{eq:modelT}).  The second system is an assembled impedance-to-impedance (ItI) trace system $\mtx{R}^{\rm ItI}$, following the construction of~\cite{2013_martinsson_ItI}.  The third system is the overlapping equilibrium system $\mtx{S}$ from (\ref{eq:modelS}).  Table~\ref{tab:example2_conditioning} reports 
unrestarted GMRES iteration counts.

As shown in Figure~\ref{fig:local_bvp}(b), let $\mtx{T}^{(\alpha)}$ and $\mtx{T}^{(\beta)}$ denote the local outward-normal DtN maps on the left and right subdomains, respectively, with $\mtx{u}$ and $\mtx{g}$ denoting the Dirichlet and Neumann data.  For impedance parameter $\eta$, define on each subdomain the incoming and outgoing Robin traces by
\[
  \mtx{r}^{(\alpha)}=\mtx{g}^{(\alpha)}-i\eta\,\mtx{u}^{(\alpha)}, \qquad
  \mtx{s}^{(\alpha)}=\mtx{g}^{(\alpha)}+i\eta\,\mtx{u}^{(\alpha)}.
\]
For the Helmholtz experiments in the code we take $\eta=\kappa$.
The local ItI maps are then defined by the Cayley transform
\[
  \mtx{R}^{(\alpha)}
  =\left(\mtx{T}^{(\alpha)}+i\eta \mtx{I}\right)
   \left(\mtx{T}^{(\alpha)}- i\eta \mtx{I}\right)^{-1},
\]
which map incoming to outgoing Robin traces, as in \cite[eqs.~17--21]{2013_martinsson_ItI}. 
The corresponding traces and ItI map on $\Omega_j^{(\beta)}$ are defined analogously.
For two neighboring subdomains sharing an interface $\Gamma_j$, if $\mtx{r}_{j}^{(\alpha)}$ and $\mtx{s}_{j}^{(\alpha)}$ are the incoming and outgoing traces on the left subdomain, and $\mtx{r}_{j}^{(\beta)}$ and $\mtx{s}_{j}^{(\beta)}$ are the corresponding traces on the right subdomain, then the shared interface relations are
\[
  \mtx{r}_{j}^{(\alpha)} = -\mtx{s}_{j}^{(\beta)}, \qquad
  \mtx{s}_{j}^{(\alpha)} = -\mtx{r}_{j}^{(\beta)}.
\]
These signs reflect the opposite outward normals on the two sides of $\Gamma_j$.  On the shared interface block, the left and right maps read
\[
  \mtx{s}_{j}^{(\alpha)}
  = \mtx{R}_{j,j-1}^{(\alpha)} \mtx{r}_{j-1}^{(\alpha)}
    + \mtx{R}_{j,j}^{(\alpha)} \mtx{r}_{j}^{(\alpha)},
  \qquad
  \mtx{s}_{j}^{(\beta)}
  = \mtx{R}_{j,j}^{(\beta)} \mtx{r}_{j}^{(\beta)}
    + \mtx{R}_{j,j+1}^{(\beta)} \mtx{r}_{j+1}^{(\beta)}.
\]
Combining these with the compatibility relations gives
\[
  \mtx{r}_{j}^{(\alpha)}
  =
  \mtx{W}\Bigl(
  \mtx{R}_{j,j}^{(\beta)}\mtx{R}_{j,j-1}^{(\alpha)}\mtx{r}_{j-1}^{(\alpha)}
  -\mtx{R}_{j,j+1}^{(\beta)}\mtx{r}_{j+1}^{(\beta)}
  \Bigr),
\qquad
  \mtx{W}
  =
  \bigl(\mtx{I}-\mtx{R}_{j,j}^{(\beta)}\mtx{R}_{j,j}^{(\alpha)}\bigr)^{-1}.
\]
This is the dense merge solve in \cite[eq.~29]{2013_martinsson_ItI} written in the left-right notation used here. 
Applying this merge at each shared slab interface eliminates the duplicate
left and right Robin traces.  We then keep one Robin trace unknown, denoted
by $\mtx{r}_j$, on each physical interface $\Gamma_j$.  The resulting assembled
ItI interface system is block tridiagonal; we denote it by $\mtx{R}^{\rm ItI}$
and use it as the ItI system in the comparison below.

Over the refinement range in Table~\ref{tab:example2_conditioning}, the DtN
Schur-complement iteration count increases from 762 to 1497.  The ItI count
increases more mildly, from 306 to 358.  The overlapping system remains at
191 GMRES iterations throughout.  Thus, for this high-frequency waveguide
example and the physical right-hand side associated with (\ref{eq:waveGuide}),
the overlapping formulation is least sensitive to mesh refinement among the
three trace systems.

\begin{table}[H]
\centering
\caption{Unrestarted GMRES counts (relative tolerance $10^{-8}$) for the three trace systems in Example~2 under mesh refinement. The slab geometry is fixed at $N_{\rm ds}=7$, and the local HPS order is fixed at $p=8$.}
\label{tab:example2_conditioning}
\begin{tabular}{c|ccc}
$N$ & $\mtx{T}$ iters & $\mtx{R}^{\rm ItI}$ iters & $\mtx{S}$ iters\\ 
\hline
589{,}824   & 762  & 306 & 191\\ 
1{,}327{,}104 & 926  & 326 & 191\\
2{,}359{,}296 & 1065 & 333 & 191 \\
9{,}437{,}184 & 1497 & 358 & 191 \\ 
\end{tabular}
\end{table}
\clearpage
\section{Spectrum of the collocated equilibrium operator}\label{app:spectra}
Here we inspect the eigenvalues of $\mtx{S}$, where $\mtx{S}$ is a sufficiently fine spectral discretization of the equilibrium operator $\mathcal{S}$, obtained from various partial differential operators.

In this section and Section~\ref{app:spectraHtozero} we follow the principle of $L^2$-weighted discretizations for continuous operators (see~\cite[\S43]{embree}): we discretize $\mathcal{S}$ to $\mtx{S}$ in such a way that $\langle\mathcal{S}u,v\rangle_{L^2(\Gamma)}\approx \vct{v}^*\mtx{S}\vct{u}$ where $\vct{u}$ and $\vct{v}$ are the discretizations of the interface functions $u$ and $v$ respectively and $\Gamma=\bigcup_j\Gamma_j$. For instance, if $\widetilde{\mtx{S}}$ denotes the (unweighted) spectrally collocated equilibrium operator, $\mtx{S}$ is obtained by scaling $\widetilde{\mtx{S}}$ by a diagonal matrix $\mtx{W} = \left(\mtx{I}_{N_{\text{ds}}}\otimes \mtx{D}_{\vct{w}}\right)$, such that $\mtx{S} = \mtx{W}\widetilde{\mtx{S}}\mtx{W}^{-1}$ with $\mtx{D}_{\vct{w}}:=\text{diag}(\vct{w})$ and $\vct{w}$ containing the square roots of the Clenshaw-Curtis quadrature weights. Any weighted discretization $\mtx{S}$ of $\mathcal{S}$ can be used to solve the original unweighted system $\widetilde{\mtx{S}}\vct{u} = \vct{b}$ by solving $\mtx{S}\vct{y} = \mtx{W}\vct{b}$ and $\vct{u} = \mtx{W}^{-1}\vct{y}$. In Figures~\ref{fig:eigValsLaplHelm} and~\ref{fig:eigValsVC} we plot the eigenvalues of $\mtx{S}$ (after correct $L^2$-weighting).

For the Laplace operator and the Helmholtz operator (at $\kappa = 10$) we see in Figure~\ref{fig:eigValsLaplHelm} that both spectra are purely real. We clearly see that for the Laplace operator the matrix $\mtx{S}$ is symmetric positive definite, and for the Helmholtz equation it is self-adjoint. We see that the Fredholm character of $\mathcal{S}=(\mathcal{I}-\mathcal{K})$ is preserved by the discretization, i.e. $\mtx{S}=\mtx{I}-\mtx{K}$. Indeed, by construction, $\mtx{K}$ acts like a discretized Hilbert-Schmidt kernel integral operator. For the Helmholtz equation, the presence of a damping term perturbs the spectrum into the negative half-plane and past $2$, which means our analysis in Section~\ref{sec:conditioning} needs to be refined in this case.
\begin{figure}[H]
    \centering
    \begin{subfigure}{.48\linewidth}
    \begin{tikzpicture}
    \begin{axis}[
    width=\linewidth,
    xlabel={$\Re (\lambda)$},
    ylabel={$\Im (\lambda)$},
    xmin=0, xmax=2,
    grid style=dashed,
    legend pos = south east
]
\addplot[only marks,mark options={.,black}
    ] table [x=real,y=imag,col sep=comma]{dataOMS/eig_Laplspectral.csv};
    \end{axis}
    \end{tikzpicture}
    \caption{Spectrum of the discretized equilibrium operator $\mtx{S}$ for the Laplace equation.}
    \end{subfigure}
\hfill
\begin{subfigure}{.48\linewidth}
    \begin{tikzpicture}
    \begin{axis}[
    /pgf/number format/.cd,fixed,precision=6,
    width=\linewidth,
    xlabel={$\Re (\lambda)$},
    ylabel={$\Im (\lambda)$},
    xmin=-1.5, xmax=3.5,
    grid style=dashed
]
\addplot[only marks,mark options={.,black}
    ] table [x=real,y=imag,col sep=comma]{dataOMS/eig_Helmspectral.csv};
    \end{axis}
    \end{tikzpicture}    
    \caption{Spectrum of the discretized equilibrium operator $\mtx{S}$ for the Helmholtz equation.}
    \end{subfigure}    
    \caption{Spectrum of $\mtx{S}$ for the Laplace equation (left) and the Helmholtz equation (right).}
    \label{fig:eigValsLaplHelm}
\end{figure}

\subsection{The variable-coefficient case}\label{app:spectraHtozero}
In Section~\ref{sec:conditioning} and the previous section, we essentially rely on the condition number estimate $\kappa_{\rho}:=\rho(\mtx{S})\rho(\mtx{S}^{-1})$ and show that $\kappa_{\rho}=\mathcal{O}(1/H^2)$. In this section we present an analysis of the continuum operator $\mathcal{S}$ underlying $\mtx{S}$ and use this to investigate $\kappa_{2}(\mathcal{S})$. Our analysis, and the supporting numerical results suggest that $\kappa_{2}(\mtx{S}) = \mathcal{O}(\kappa_{\rho})$. Intuitively, this means that the measure of non-normality of the matrix $\mtx{S}$ is benign.

The outline of our argument is as follows:
\begin{enumerate}
    \item We assume that we have a square domain $\Omega=[0,1]^2$, and equispaced interfaces $\{\Gamma_j\}_{j=1}^{N_{\text{ds}}}$ with separation $H$.
    \item We show that, as $H\rightarrow 0$, the continuous equilibrium operator $\mathcal{S}$ becomes in a sense `sufficiently self-adjoint'.
    \item Then we have that, for an $L^2$-weighted discretization $\mtx{S}$ of $\mathcal{S}$, this implies that $\mtx{S}$ also becomes sufficiently self-adjoint.
\end{enumerate}

For ease of introduction we first assume that we have a positive definite formally self-adjoint PDE operator $\mathcal{A}$ on $\Omega$ with \textit{constant coefficients}. Our first task is then to show that in this case the equilibrium operator $\mathcal{S}$ is self-adjoint. It is clear from the continuum form in \ref{eq:debbie1} that
$$[\mathcal{S}^*_{j-1,j}u_{j-1}](x)=\int_{\Gamma_{j-1}}\overline{G^{(j-1)}(y,x)}u_{j-1}(y)dy=\int_{\Gamma_{j-1}}G^{(j)}(x,y)u_{j-1}(y)dy$$
since $G^{(j-1)}(x,y)=\overline{G^{(j-1)}(y,x)}$ and $G^{(j)}=G^{(j-1)}$. This last claim is true only by virtue of $\mathcal{A}$ having constant coefficients \textbf{and} the fact that all slabs are chosen isomorphic. Therefore, in the simple case considered, we have that $\mathcal{S} = \mathcal{S}^*$.

For the case of non-constant coefficients the operator $\mathcal{S}$ is no longer self-adjoint. It is not even a normal operator. However, as $H\rightarrow 0$ we can still recover that $G^{(j)}\rightarrow G^{(j-1)}$ and vice versa. The rate at which this happens depends on the smoothness of the coefficients of $\mathcal{A}$. As such, for formally self-adjoint PDE operators with sufficiently smooth coefficients, an \textit{asymptotic version} of self-adjointness can still be recovered. While a full analysis is beyond the scope of this work, we mention that this is essentially because the Green's function in the case of smooth coefficients varies smoothly with the \textit{perturbation} of moving from $\Psi_j$ to $\Psi_{j-1}$. This can be shown using the techniques from \cite{GarabedianSchiffer}, \S \RN{2} (see also \cite{Kato}, \S \RN{7}.6.5). We will present a numerical study of the asymptotics here.

We use the variable-coefficient convection-diffusion differential operator from~\eqref{eq:screened_diff} in Section~\ref{sec:convdiff} with $\alpha = 1$, but let us mention that the same behavior was observed for any other positive definite elliptic variable-coefficient PDE operator tested. For a given $H$, we construct two approximations $\mtx{S}_{H,h}$ and $\mtx{S}_{H,p}$ of the corresponding $\mathcal{S}$; respectively these are built with a fine five-point stencil discretizations (isotropic, with $h_x=h_y=h=1/128$) and with high-order Chebyshev spectral collocation discretizations (with order $p_x=p_y=p=30$) for the overlapping slabs. As described above, $\mtx{S}_{H,p}$ is $L^2$-weighted. The discretization $\mtx{S}_{H,h}$ does not need to be weighted.

In Figure~\ref{fig:normalFigAsymptotic}, we investigate for both discretizations three measures of normality. Firstly, we compute $\|\mtx{S}_{j,j-1}-\mtx{S}_{j-1,j}^*\|_2$ with $\Gamma_j$ a fixed interface, as an indicator of the smoothness of the Green's functions over $H$. In our case we chose $\Gamma_{j}$ to correspond to the interface $x=1/2$, meaning $\Gamma_{j-1}$ corresponds to $x=1/2-H$. Secondly, we compute
$$\|\vct{\lambda}-\vct{\sigma}\|_{\infty}:=\max_i\{{||\lambda_i|-\sigma_i}|\}$$
where $|\lambda_1|,|\lambda_2|,\ldots$ and $\sigma_1\geq \sigma_2\geq \cdots$ are the moduli of the eigenvalues and the singular values of the discretizations\footnote{The ordering of the eigenvalues is chosen so as to minimize $\|\vct{\lambda}-\vct{\sigma}\|_{\infty}$.}. Finally, we also plot the measure $\kappa_2/\kappa_{\rho}-1$, where $\kappa_{\rho}:=\rho(\mtx{S}_{H,h})\rho(\mtx{S}_{H,h}^{-1})$ (similarly for $\mtx{S}_{H,p}$), and $\kappa_2$ is the $\|\cdot\|_2$-condition number.

\begin{figure}[H]
    \pgfplotsset{select coords between index/.style 2 args={
    x filter/.code={
        \ifnum\coordindex<#1\def\pgfmathresult{}\fi
        \ifnum\coordindex>#2\def\pgfmathresult{}\fi
    }
}}
\centering
    \begin{subfigure}[t]{.32\linewidth}
    \centering
    \begin{tikzpicture}[every node/.append style={font=\small}]
    \begin{axis}[
    width=\linewidth,
    height=6cm,
    xmode=log,
    ymode=log,
    log basis x=2,
    log basis y=10,
    xlabel={$H$},
    ylabel={\small$\|\mtx{S}_{j,j-1}-\mtx{S}_{j-1,j}^*\|_{2}$},
    ymajorgrids=true,
    xmajorgrids=true,
    grid style=dashed,
    legend pos = south east
]
\addplot[
    color=black,
    mark=triangle*,
    ] table [x=H,y=nrm_diff_stencil,col sep=comma]{dataOMS/cond_stats_conv_diff.csv};
    \addlegendentry{$\mtx{S}_{H,h}$}
\addplot[
    color=black,
    mark=square,
    ] table [x=H,y=nrm_diff_HPS,col sep=comma]{dataOMS/cond_stats_conv_diff.csv};
    \addlegendentry{$\mtx{S}_{H,p}$}
\addplot[
    color=black,
    dashed
    ] table [select coords between index={1}{5},x=H,y expr=(x),col sep=comma]{dataOMS/cond_stats_conv_diff.csv}node[pos=.25] (OH) {};
    \node [left] at (OH) {$\mathcal{O}(H)$};
    \end{axis}
    \end{tikzpicture}
    \caption{$\|\mtx{S}_{j,j-1}-\mtx{S}_{j-1,j}^*\|_{2}$ as a function of $H$.}
    \label{fig:blockErr}
    \end{subfigure}
    \hfill
    \begin{subfigure}[t]{.32\linewidth}
    \centering
    \begin{tikzpicture}[every node/.append style={font=\small}]
    \begin{axis}[
    width=\linewidth,
    height=6cm,
    xmode=log,
    ymode=log,
    log basis x=2,
    log basis y=10,
    xlabel={$H$},
    ylabel={\small$\|\vct{\lambda}-\vct{\sigma}\|_{\infty}$},
    xmin=0, xmax=.5,
    ymajorgrids=true,
    xmajorgrids=true,
    grid style=dashed,
    legend pos = south east
]
\addplot[
    color=black,
    mark=triangle*,
    ] table [x=H,y=lamdiff_stencil,col sep=comma]{dataOMS/cond_stats_conv_diff.csv};
    \addlegendentry{$\mtx{S}_{H,h}$}
\addplot[
    color=black,
    mark=square,
    ] table [x=H,y=lamdiff_HPS,col sep=comma]{dataOMS/cond_stats_conv_diff.csv};
    \addlegendentry{$\mtx{S}_{H,p}$}
\addplot[
    color=black,
    dashed
    ] table [select coords between index={1}{5},x=H,y expr=(x^1.5)/16,col sep=comma]{dataOMS/cond_stats_conv_diff.csv}node[pos=.25] (OH) {};
    \node [left] at (OH) {$\mathcal{O}(H^{3/2})$};
    \end{axis}
    \end{tikzpicture}
    \caption{$\|\vct{\lambda}-\vct{\sigma}\|_{\infty}$ as a function of $H$.}
    \label{fig:eigSvdH}
    \end{subfigure}
    \hfill
\begin{subfigure}[t]{.32\linewidth}
    \centering
    \begin{tikzpicture}[every node/.append style={font=\small}]
    \begin{axis}[
    /pgf/number format/.cd,fixed,precision=6,
    width=\linewidth,
    height=6cm,
    xmode=log,
    log basis x=2,
    xlabel={$H$},
    ylabel={\small$\kappa_2/\kappa_{\rho}-1$},
    xmin=0, xmax=.5,
    ymajorgrids=true,
    xmajorgrids=true,
    grid style=dashed,
    legend pos = south west
]
\addplot[
    color=black,
    mark=triangle*,
    ] table [x=H,y =condstat_stencil,col sep=comma]{dataOMS/cond_stats_conv_diff.csv};
    \addlegendentry{$\mtx{S}_{H,h}$}
\addplot[
    color=black,
    mark=square,
    ] table [x=H,y =condstat_HPS,col sep=comma]{dataOMS/cond_stats_conv_diff.csv};
    \addlegendentry{$\mtx{S}_{H,p}$}
    \end{axis}
    \end{tikzpicture}
    \caption{$\kappa_2/\kappa_{\rho}-1$ as a function of $H$.}
    \label{fig:kappaSvdEig}
    \end{subfigure}  
\caption{The normality measures $\|\mtx{S}_{j,j-1}-\mtx{S}_{j-1,j}^*\|_{2}$, $\|\vct{\lambda}-\vct{\sigma}\|_{\infty}$ and $\kappa_{\sigma}/\kappa_\rho-1$ as a function of the slab width $H$, for both a stencil and spectral discretization.}
\label{fig:normalFigAsymptotic}
\end{figure}
Figure~\ref{fig:normalFigAsymptotic} shows that both discretizations behave completely similarly. This is not surprising, since at each $H$ they approximate the same continuous Fredholm type operator. We see that the measure of non-normality $\|\vct{\lambda}-\vct{\sigma}\|_{\infty}$ approaches zero as $H\rightarrow 0$. In fact we have the stronger observation that $\|\vct{\lambda}-\vct{\sigma}\|_{\infty}=\mathcal{O}(H^{3/2})$. This does not show, however, that $\kappa_{\rho}\rightarrow \kappa_2$, since $\kappa_\rho$ and $\kappa_2$ both diverge as $H$ approaches zero. However, Figure~\ref{fig:kappaSvdEig} shows that $\kappa_2=c_H\kappa_{\rho}$ for the discretizations tested here, with $c_H$ remarkably close to one. Let us mention again that we observed the same behavior for every symmetric positive elliptic PDE discretization we tested.

The experiments and analysis here therefore support the use of $\kappa_{\rho}$ as an indicator of the actual condition number for the stencil and spectral discretizations tested here, and suggest similar asymptotic behavior for the two discretizations.

In Figure~\ref{fig:eigValsVC} we plot the eigenvalues of $\mtx{S}_{H,p}$ at $H=1/16$. It shows that $\mtx{S}_{H,p}$ is no-longer self-adjoint; its spectrum has a small but non-negligible imaginary component. Note that the eigenvalues not only occur in conjugate pairs (as is to be expected), but they are also symmetric around $\Re(z)=1$.

\begin{figure}[H]
    \centering
    \begin{subfigure}{.48\linewidth}
    \begin{tikzpicture}
    \begin{axis}[
    width=\linewidth,
    xlabel={$\Re (\lambda)$},
    ylabel={$\Im (\lambda)$},
    xmin=0, xmax=2,
    grid style=dashed,
    legend pos = south east
]
\addplot[only marks,mark options={.,black}
    ] table [x=real,y=imag,col sep=comma]{dataOMS/eig_spectral.csv};
    \end{axis}
    \end{tikzpicture}    
    \end{subfigure}
\hfill
\begin{subfigure}{.48\linewidth}
    \begin{tikzpicture}
    \begin{axis}[
    /pgf/number format/.cd,fixed,precision=6,
    width=\linewidth,
    xlabel={$\Re (\lambda)$},
    ylabel={$\Im (\lambda)$},
    xmin=.9995,xmax=1.0005,
    xtick={.9995,1,1.0005},
    grid style=dashed
]
\addplot[only marks,mark options={.,black}
    ] table [x=real,y=imag,col sep=comma]{dataOMS/eig_spectral_im.csv};
    \end{axis}
    \end{tikzpicture}    
    \end{subfigure}    
    \caption{Total spectrum (left) and non-real part of the spectrum (right) of the discretized equilibrium operator $\mtx{S}$ for a variable-coefficient convection-diffusion equation.}
    \label{fig:eigValsVC}
\end{figure}
\subsection{Comparison to the non-overlapping domain decomposition}\label{app:Tsys}
Here we investigate the spectrum of the $\mtx{T}$ system. For compactness we restrict our attention to the Laplace equation. We use fine stencil and spectral discretizations and report the results in Figure~\ref{fig:spectrumT}. In both cases the domain $\Omega=[0,1]^2$ is the unit square, which was divided into non-overlapping slabs of width $H=1/16$. For the spectral case we employ the $L^2$-weighting principle outlined in Section~\ref{app:spectra}.

We clearly see that the spectral range is not only much larger than that of the $\mtx{S}$-system, but also that the spectral behavior of the underlying continuum operator $\mathcal{T}$ is not captured well by both discretizations. In short, the operator $\mathcal{T}$ behaves like an unbounded pseudo-differential operator, as opposed to the second kind Fredholm behavior of the operator $\mathcal{S}$. For the stencil, we see that the largest eigenvalue grows like $\mathcal{O}(1/h^2)$, where similarly the largest eigenvalue for the spectral case scales like $\mathcal{O}(p^2)$.
\begin{figure}[H]
    \centering
    \begin{subfigure}{.48\linewidth}
    \begin{tikzpicture}
    \begin{axis}[
    width=\linewidth,
    height=4cm,
    xlabel={$\Re (\lambda)$},
    ylabel={$\Im (\lambda)$},
    grid style=dashed,
    legend pos = south east
]
\addplot[only marks,mark options={.,black}
    ] table [x=real,y=imag,col sep=comma]{dataOMS/spectrumT_0.030303030303030304.csv};
    \end{axis}
    \end{tikzpicture}
    \begin{tikzpicture}
    \begin{axis}[
    width=\linewidth,
    height=4cm,
    xlabel={$\Re (\lambda)$},
    ylabel={$\Im (\lambda)$},
    grid style=dashed,
    legend pos = south east
]
\addplot[only marks,mark options={.,black}
    ] table [x=real,y=imag,col sep=comma]{dataOMS/spectrumT_0.016.csv};
    \end{axis}
    \end{tikzpicture}
    \begin{tikzpicture}
    \begin{axis}[
    width=\linewidth,
    height=4cm,
    xlabel={$\Re (\lambda)$},
    ylabel={$\Im (\lambda)$},
    grid style=dashed,
    legend pos = south east
]
\addplot[only marks,mark options={.,black}
    ] table [x=real,y=imag,col sep=comma]{dataOMS/spectrumT_0.007751937984496124.csv};
    \end{axis}
    \end{tikzpicture}
    \caption{Spectrum of the matrix $\mtx{T}$ for the Laplace equation with stencil discretizations at $h=2^{-k}$, $k=6$ (top), $k=7$ (middle) and $k=8$ (bottom). Note the difference in length scales for the real part of the spectrum.}
    \end{subfigure}
\hfill
    \begin{subfigure}{.48\linewidth}
    \begin{tikzpicture}
    \begin{axis}[
    width=\linewidth,
    height=4cm,
    xlabel={$\Re (\lambda)$},
    ylabel={$\Im (\lambda)$},
    grid style=dashed,
    legend pos = south east
]
\addplot[only marks,mark options={.,black}
    ] table [x=real,y=imag,col sep=comma]{dataOMS/spectrumT_spectral64.csv};
    \end{axis}
    \end{tikzpicture}
    \begin{tikzpicture}
    \begin{axis}[
    width=\linewidth,
    height=4cm,
    xlabel={$\Re (\lambda)$},
    ylabel={$\Im (\lambda)$},
    grid style=dashed,
    legend pos = south east
]
\addplot[only marks,mark options={.,black}
    ] table [x=real,y=imag,col sep=comma]{dataOMS/spectrumT_spectral128.csv};
    \end{axis}
    \end{tikzpicture}
    \begin{tikzpicture}
    \begin{axis}[
    width=\linewidth,
    height=4cm,
    xlabel={$\Re (\lambda)$},
    ylabel={$\Im (\lambda)$},
    grid style=dashed,
    legend pos = south east
]
\addplot[only marks,mark options={.,black}
    ] table [x=real,y=imag,col sep=comma]{dataOMS/spectrumT_spectral256.csv};
    \end{axis}
    \end{tikzpicture}
    \caption{Spectrum of the matrix $\mtx{T}$ for the Laplace equation with spectral discretizations at $p=64$ (top), $p=128$ (middle) and $p=256$ (bottom). Note the difference in length scales for the real part of the spectrum.}
    \end{subfigure}
    \caption{Spectrum of the matrix $\mtx{T}$ for the Laplace equation at different discretizations (stencil and spectral). Observe that the spectrum (its bounds \textbf{and} its shape) is discretization dependent. Additionally, there is no beneficial spectral clustering as in the case for the $\mtx{S}$-matrix.}
    \label{fig:spectrumT}
\end{figure}

\end{document}